\documentclass[11pt,a4paper]{article}
\usepackage[T1]{fontenc}
\usepackage[utf8]{inputenc}
\usepackage{lmodern}
\usepackage[a4paper,margin=0.70in]{geometry}
\usepackage{amsmath,amssymb,mathrsfs,bm,mathtools}
\usepackage{graphicx,booktabs,float,subfig,multirow,tabularx,array}
\usepackage{enumitem,tikz,changepage}
\usepackage{microtype}
\usepackage[numbers,sort&compress]{natbib}
\usepackage[hidelinks]{hyperref}
\usepackage{authblk}
\usepackage{xcolor}
\usetikzlibrary{shapes.geometric,arrows.meta,positioning,arrows,decorations.pathreplacing,patterns}
\numberwithin{equation}{section}
\allowdisplaybreaks
\newcolumntype{C}{>{\centering\arraybackslash}X}
\newcolumntype{L}{>{\raggedright\arraybackslash}X}
\newlength{\extralength}
\newcommand{\fulllength}{\textwidth}
\newcommand{\cbl}[1]{#1}
\newcommand{\linksupplementary}[1]{\href{https://doi.org/10.3390/fractalfract10110745}{the published article page}}

\title{Equation-Free Screening of Mittag--Leffler Compatibility\\from Multi-Horizon Relaxation Profiles}
\author[1,2]{N'Gbo N'Gbo}
\author[3]{Andrei Velichko}
\author[4]{Bruno Carpentieri\thanks{Correspondence: \href{mailto:bruno.carpentieri@unibz.it}{bruno.carpentieri@unibz.it}}}
\affil[1]{School of Science and Engineering, International University of Grand-Bassam, Route de Bonoua, Grand-Bassam BP 564, C\^ote d'Ivoire}
\affil[2]{Laboratory for Intelligence and Mathematics (LIMAs), Route de Bonoua, Grand-Bassam BP 564, C\^ote d'Ivoire}
\affil[3]{Institute of Physics and Technology, Petrozavodsk State University, 33 Lenina Prospekt, Petrozavodsk 185910, Republic of Karelia, Russia}
\affil[4]{Faculty of Engineering, Free University of Bozen-Bolzano, Via Bruno Buozzi 1, 39100 Bozen--Bolzano, Italy}
\date{26 October 2026}

\hypersetup{
  pdftitle={Equation-Free Screening of Mittag-Leffler Compatibility from Multi-Horizon Relaxation Profiles},
  pdfauthor={N'Gbo N'Gbo, Andrei Velichko, Bruno Carpentieri}
}

\begin{document}
\maketitle
\begin{center}
\footnotesize
N'Gbo N'Gbo: ORCID \href{https://orcid.org/0009-0005-3056-4688}{0009-0005-3056-4688} \quad
Andrei Velichko: ORCID \href{https://orcid.org/0000-0002-9341-1831}{0000-0002-9341-1831}\\
Bruno Carpentieri: ORCID \href{https://orcid.org/0000-0002-0516-0033}{0000-0002-0516-0033}
\end{center}

\begin{abstract}

We study an equation-free screen for finite-window Mittag--Leffler compatibility,
	in which multi-horizon relaxation profiles are constructed out of sample by
	k-nearest-neighbor prediction from ensembles of scalar trajectories and compared with
	candidate temporal laws. The screen assesses compatibility, not the generating mechanism.
	In a matched, noise-free calibration, it reproduces the generating Mittag--Leffler behavior
	to numerical precision, while relaxation-rate heterogeneity shifts the fitted
	Mittag--Leffler order toward lower values. For noise-free classical distributed-rate
	ensembles, a Mittag--Leffler preference over a single exponential disappears when a
	biexponential alternative is included. Under observational noise, this recovery is
	incomplete: with three candidate families, the Mittag--Leffler law is still selected in
	60\% of runs for a classical ensemble, against 72\% for a fractional reference. In the
	tested noise-free settings, simpler profile constructions reproduce the screening
	decisions, whereas kNN predicts individual trajectories far more accurately than a
	constant predictor. The screening interpretation of a single-record example, initially
	investigated as a chaotic test case, is withdrawn as evidence about the continuous system
	after numerical verification revealed a nonrobust simulated regime. Mittag--Leffler
	compatibility is therefore a comparative, condition-dependent screening outcome, not
	identification of fractional origin.
\end{abstract}


\noindent\begin{minipage}{\textwidth}\footnotesize\itshape
Author-prepared version corresponding to the published article:\par
\emph{Fractal and Fractional} 2026, 10(11), 745.\quad
\url{https://doi.org/10.3390/fractalfract10110745}\par
Reproducibility package: \url{https://doi.org/10.5281/zenodo.23494805}.
\end{minipage}\par
\normalsize\medskip

\section{Introduction}
\label{Sec:Introduction}

Fractional calculus provides a natural mathematical framework for systems in which the
present state depends on a distributed history rather than only on an instantaneous local
state. Caputo and related fractional operators have therefore been used to represent
hereditary effects, anomalous relaxation, viscoelasticity, transport, electrochemical memory,
and nonlinear dynamical behavior
~\cite{podlubny1999fractional,diethelm2010analysis,Herrmann2014Fractional,Mainardi2000On,Mainardi2020Why}.
A canonical case is the linear Caputo relaxation problem
${}_{C}D_{0,t}^{\alpha}x=-\lambda x$, whose solution is the Mittag--Leffler law
$x(t)=x_0E_\alpha(-\lambda t^\alpha)$, with~the ordinary exponential recovered in the
integer-order limit $\alpha=1$
\cite{Mainardi2000On,Mainardi2020Why,Gorenflo2020MittagLeffler}.
In this canonical setting, the~governing fractional dynamics are linked analytically to a
scalar relaxation~response.

This explicit link suggests a data-driven question. When the governing equations are
unknown and only scalar observations are available, can the data reveal a temporal pattern
compatible with this canonical fractional response? Such compatibility could motivate
further investigation of fractional~models.

The difficulty is that such a pattern is not unique to fractional dynamics. Relaxation
curves that resemble a Mittag--Leffler law can also be produced by classical mechanisms: stretched exponentials,
mixtures or distributions of ordinary exponential rates, aggregation, structural changes,
and nonstationarity can all generate similarly broad or curved responses over a finite
observation window~\cite{Metzler2002From,Granger1980AN,Monache2015Testing,Mainardi2017A}, and~observational
noise can further obscure their differences.
Agreement with a Mittag--Leffler-type pattern can therefore indicate compatibility with a
fractional description, but~cannot by itself identify the mechanism that generated the
data. It is useful to separate
three questions that are often merged. The~first is whether the observed temporal
behavior is compatible with a fractional-type law. The~second is whether a fractional generating mechanism
is better supported than plausible classical alternatives. The~third
is what the parameters of that mechanism are, including the fractional order. Each question
demands stronger evidence than the one before~it.

Existing fractional-system identification methods work within a specified fractional
model class: they estimate its order, coefficients, hidden states, or~transfer function
from input--output or state data
~\cite{Poinot2004Identification,Li2021Fractional,Zhang2024Sparse,
	duhe2024recursive,Yuan2024Identification}.
Related model-based work also addresses observability and initial-state reconstruction
in time-fractional systems~\cite{BenBrahim2024Some}.
These methods are powerful when a credible governing model is available, but~they operate
after a fractional model class has already been posited. In~many applications only a single
scalar quantity is measured and no governing ODE, PDE, or~FDE has been selected. In~such
settings, a~preliminary data-driven assessment of the first question can be useful before
committing to a particular governing model
~\cite{Du2013Measuring,Tarasov2019On,Wei2024Identifying}. The~purpose of this work is to
provide such an assessment: an earlier, equation-free screening layer for settings in which
the appropriate governing model class is still~uncertain.

The screening principle builds on the use of multi-horizon prediction as an equation-free
probe of temporal dynamics, an~idea we previously exploited to estimate largest Lyapunov
exponents from one-dimensional chaotic time series~\cite{Velichko2025A}. Rather than
fitting a governing differential equation to the observations, we use direct-horizon kNN
prediction to construct, out of sample, an~empirical relaxation profile from an ensemble of
scalar trajectories: each test trajectory is represented only by its observed history, and~its future is predicted from training trajectories with similar histories. Candidate
temporal laws are then fitted to this profile as reference families to quantify
finite-window compatibility, not to identify the governing equations directly. If~the
Mittag--Leffler law describes the profile better than the classical alternatives
considered, a~fractional description may warrant further model-based investigation.
Because screening is performed on the aggregate profile, it is not evident how much
the history-conditioned kNN construction contributes to the screening decision compared
with simpler alternatives, including a constant training-mean predictor and the directly
observed ensemble average; we examine this question~explicitly.

The practical value of such a screen depends on understanding both its reliability and its
limitations. The~aim of this study is therefore twofold: to propose the screen and to
examine its reliability and limitations under increasingly demanding conditions. We
address four research~questions:
\begin{enumerate}[label=,leftmargin=1cm,labelsep=3mm]
	\item[RQ1] \emph{Calibration.} Can a multi-horizon kNN profile recover the expected
	Mittag--Leffler compatibility in a controlled matched-family relaxation setting?
	\item[RQ2] \emph{Specificity and candidate-set dependence.} How do relaxation-rate
	heterogeneity, nonfractional curved mechanisms, and~observational noise affect an
	apparent Mittag--Leffler preference, and~how does that preference change when the
	candidate set \mbox{is broadened}?
	\item[RQ3] \emph{Interpretation and robustness.} How robust are the screening decision
	and the fitted Mittag--Leffler shape parameter to observational noise and analysis
	choices? When the data do not arise from a single matched law, how should the fitted
	parameter \mbox{be interpreted}?
	\item[RQ4] \emph{Profile construction.} To what extent does the kNN construction change
	the screening decision relative to simpler profile constructions, and~what does it add
	at the level of individual-trajectory prediction?
\end{enumerate}

As a controlled foundation, we deliberately begin with the canonical linear Caputo
relaxation problem introduced above. Its simplicity is precisely what makes it useful: the
generating law is known analytically, so the screening principle can first be calibrated
under matched conditions and then challenged progressively by rate heterogeneity,
nonfractional alternatives, observational noise, and~sensitivity to analysis choices, while
the kNN construction is compared with simpler profile constructions. We also report a
single-record example generated numerically from a nonlinear fractional system initially
investigated as a chaotic test case. Verification with an independent predictor--corrector
implementation showed that the sustained oscillatory regime produced by the original
explicit discretization is not numerically robust; we therefore retain this example only
as a cautionary illustration of how numerical artifacts in simulated observations can
affect data-driven~screening.

Together, these experiments probe the conditions under which a finite-window compatibility
screen is informative and those under which its interpretation becomes ambiguous. The~resulting framework is
positioned between descriptive time-series analysis and full fractional-system
identification: it is intended to motivate subsequent mechanistic modeling, not to replace
it.

Section~\ref{Sec:Methods} describes the screening principle, the~benchmarks, and~the
protocols. \mbox{Section~\ref{Sec:Results}} reports the results, and~Section~\ref{Sec:Discussion} discusses their implications and
limitations. Section~\ref{Sec:Conclusions} concludes.

\section{Materials and~Methods}
\label{Sec:Methods}

Section~\ref{Sec:ScreeningPrinciple} states the screening procedure and the separation
between model fitting and evaluation (Figure~\ref{Fig:DataSeparationSchematic}).
Sections~\ref{Sec:StableBenchmarks} and~\ref{Sec:StableProfile} define the stable
benchmarks and their observable, together with the simpler profile constructions used for
comparison (RQ1, RQ2, RQ4). Section~\ref{Sec:CandidateModels} specifies the candidate
families and their optimization, and \mbox{Sections~\ref{Sec:NoiseProtocol}
and~\ref{Sec:SensitivityProtocol}} the noise, specificity, and~sensitivity analyses (RQ2,
RQ3). Section~\ref{Sec:NumericalCheck} describes the numerical-reliability check of a
single-record~example.

\begin{figure}[H]
	
\begin{adjustwidth}{0pt}{0pt}
\centering
	\begin{tikzpicture}[x=1cm,y=1cm,font=\small]
		
		\node[anchor=west,font=\small\bfseries] at (-0.2,2.6)
		{Stable experiments: fitting and held-out horizons};
		
		\draw[->] (0,0) -- (14.4,0) node[right] {$t$};
		\foreach \t/\x in {0/0, 0.6/2.33, 2.6/10.11, 3.6/14}
		\draw (\x,0.08) -- (\x,-0.08) node[below] {\footnotesize $\t$};
		
		\fill[gray!35]   (0,1.0)     rectangle (2.33,1.5);
		\fill[blue!30]   (2.72,1.0)  rectangle (10.11,1.5);
		\fill[orange!45] (10.50,1.0) rectangle (14.0,1.5);
		
		\node at (1.17,1.25)  {\scriptsize history ($m_{\mathrm{hist}}=7$)};
		\node at (6.42,1.25)  {\footnotesize fitting horizons $h=1$--$20$};
		\node at (12.25,1.25) {\footnotesize held-out $h=21$--$30$};
		
		\draw[dashed] (2.33,-0.1) -- (2.33,1.9)
		node[above] {\footnotesize $t_{\mathrm{ref}}$};
		\draw[dashed,thick] (10.31,-0.1) -- (10.31,1.9)
		node[above] {\footnotesize fit\,/\,held-out};
		
		\node[anchor=west] at (0,-0.85)
		{\footnotesize $t_h=t_{\mathrm{ref}}+h\,\Delta t$;
			trajectories split 70/30 into training and test sets};
		
	\end{tikzpicture}
\end{adjustwidth}
	
	\caption{Separation of training and test data in the stable experiments: trajectories
		are split 70/30 into training and test sets; the candidate relaxation laws are fitted
		on horizons $h=1$--$20$ and evaluated on the held-out horizons $h=21$--$30$, on~the
		absolute target time $t_h=t_{\mathrm{ref}}+h\,\Delta t$.}
	\label{Fig:DataSeparationSchematic}
\end{figure}

\subsection{Screening Principle and Data-Separation~Protocols}
\label{Sec:ScreeningPrinciple}

The screening procedure has three steps. First, an~empirical multi-horizon relaxation
profile is constructed from an ensemble of scalar trajectories by direct-horizon kNN
prediction. Second, candidate temporal laws are fitted to this profile over a prescribed
fitting interval. Third, the~fitted laws are compared on that interval and on later
horizons not used for~fitting.

Because the screen is intended to characterize out-of-sample behavior, the~profile must
be constructed without information from the evaluation data entering the predictor. The~kNN predictor and all standardization quantities therefore use training data only.
Forecasting is also \emph{direct-horizon}: each horizon \(h\) is predicted
separately from the history vector, rather than recursively from predictions at shorter
horizons. The~profile at horizon \(h\) thus reflects the empirical predictive relationship
at that horizon, not errors accumulated by iterating a one-step~predictor.

In the stable branch, the~data are separated at the trajectory level: \(70\%\) of the
trajectories train the kNN predictor, and~the remaining \(30\%\) are used to compute the
profile for horizons \(h=1,\ldots,30\). Candidate laws are fitted on horizons
\(1\)--\(20\) and then, with~their parameters held fixed, evaluated on horizons
\(21\)--\(30\). We call this second stage a \emph{horizon-extrapolation check}: it asks
whether a law fitted to the earlier part of the profile continues to describe its later
part. Because~all horizons are computed from the same test trajectories, they are
statistically dependent, and~the later horizons are not treated as an independent
validation~sample.

Within each experiment, all candidate laws are fitted to the same profile with the same
least-squares criterion (Section~\ref{Sec:CandidateModels}). Because~the screen asks which
reference family continues to describe the profile beyond the horizons used to fit it, the~primary comparison criterion for the baseline stable experiments is the
horizon-extrapolation RMSE on horizons \(21\)--\(30\). The~corresponding intervals are
adjusted in the sensitivity analyses described in
Section~\ref{Sec:SensitivityProtocol}. In-sample criteria need not agree with the
extrapolation criterion: a more flexible family can fit the fitting horizons more closely
and yet extrapolate less well. AIC and BIC computed on the fitting horizons are therefore
reported only as descriptive complexity-penalized summaries, not as formal
likelihood-based tests, because~neighboring horizons are statistically dependent.
Neither held-out RMSE nor the in-sample criteria can identify the generating mechanism:
they assess compatibility with the constructed~profile.

The screening decision is therefore comparative rather than a formal test of fractional
origin. A~lower fitting or extrapolation error indicates compatibility with a candidate
family relative to the alternatives considered; it does not identify a unique generating
mechanism or establish a physical fractional~order.

\subsection{Stable Benchmark~Design}
\label{Sec:StableBenchmarks}

To interpret what the screen detects, it must first be tested on data whose generating law
is known. The~stable benchmarks therefore start from the most favorable situation and
progressively relax or replace its idealizations. Three principal cases separate three questions: exact
matched-family calibration (Case A), robustness to moderate rate heterogeneity within a
fractional relaxation family (Case B), and~specificity against a purely classical
distributed-rate mechanism (Case C). Two additional shape controls probe a different
issue: whether strong curvature alone can generate an apparent Mittag--Leffler-compatible
profile, and~whether a correctly specified curved classical family is recovered when it is
included among the~candidates.

\subsubsection{Principal Stable Benchmark~Sequence}
\label{Sec:StableABC}

For Cases A and B, trajectories are generated from the analytical fractional
relaxation law 
\begin{equation}
	x_j(t)=x_{0,j} E_{\alpha}\!\left(-\mu_j t^{\alpha}\right),
	\label{Eq:stable_ml_generator}
\end{equation}
where \(E_{\alpha}\) denotes the one-parameter Mittag--Leffler function and
\(x_{0,j}\sim\mathcal{U}[0.5,3]\).  We use \(N_{\mathrm{traj}}=500\) trajectories for
each parameter setting.  The~Mittag--Leffler function is evaluated with the
robust algorithm of Garrappa~\cite{Garrappa2015Numerical}, so no time-integration discretization is introduced
in the stable data generation.  Only the interval \(t\in[0,3.6]\), spanned by the history window and the 30
prediction horizons (Section~\ref{Sec:StableProfile}), enters the analysis; samples beyond
\(t=3.6\) that the scripts may generate are not~used.

\paragraph{Case A:  exact matched-family calibration.}
Case A is the best-case reference. The~data are generated from the same
Mittag--Leffler family that is later fitted, with~\(\alpha\in\{0.6,0.7,0.8,0.9,1\}\) and a common relaxation coefficient
\begin{equation}
	\mu_j=\mu=0.15
\end{equation}
for all trajectories.  Because~\(\mu\) is common to all trajectories,
\[
x_j(t)=x_{0,j} f(t),
\qquad
f(t)=E_{\alpha}\!\left(-\mu t^{\alpha}\right),
\]
so every history is a rescaled copy of the same temporal template.  Under~this
separable structure, the~trajectory-dependent amplitude factor introduced by
the kNN averaging is absorbed by the fitted amplitude \(C\).  Near-exact recovery
is therefore the expected behavior here: Case A provides a calibration check of the
numerical screening pipeline under exact controlled conditions, and~its success is not evidence of
general~identifiability.

\paragraph{Case B: fractional relaxation with moderate rate heterogeneity.}
Case A is deliberately ideal.  Case B removes its multiplicative separability
while remaining close to the reference relaxation scale: each trajectory is
assigned its own rate,
\begin{equation}
	\mu_j\sim\mathcal{U}[0.05,0.25],
	\label{Eq:case_b_rate}
\end{equation}
with a common fractional order \(\alpha\).  The~interval is centered at
\(0.15\) and introduces moderate inter-trajectory heterogeneity without changing
the underlying fractional relaxation family.  Case B therefore tests whether the
Mittag--Leffler-compatible profile remains detectable when the ensemble can no
longer be written as \(x_j(t)=x_{0,j}f(t)\) with one common temporal factor
\(f\).

\paragraph{Case C: classical distributed-rate challenge.}
Case B still lies within a fractional generating family. Case C asks the stronger
specificity question: can rate heterogeneity alone produce the same apparent
signature when every individual trajectory is classical?  Each trajectory is purely
exponential,
\begin{equation}
	x_j(t)=x_{0,j}\exp(-\mu_j t),
\end{equation}
but the rates are log-uniformly distributed according to
\begin{equation}
	\mu_j=0.15\,s^{u_j},
	\qquad
	u_j\sim\mathcal{U}[-1,1],
	\qquad
	s\in\{3,5,10\}.
	\label{Eq:case_c_rate}
\end{equation}
The distribution is centered geometrically at the same reference value
\(0.15\), while \(s\) controls the spread.  The~spread is deliberately much wider
than in Case B: Case C is not a perturbation of the reference configuration but an
increasingly demanding classical control as \(s\) grows.  Case C is first
used in the same binary exponential-versus-Mittag--Leffler comparison as Cases A
and B and is then revisited with a broader candidate set that includes a
biexponential classical alternative; the fitting families themselves are defined
in Section~\ref{Sec:CandidateModels}.

\subsubsection{Separable Nonfractional Shape~Controls}
\label{Sec:StableShapeControls}

Cases A--C vary the rate structure of the ensemble.  The~two shape controls ask a
different question: can a nonfractional but strongly curved temporal law mimic
Mittag--Leffler compatibility, and~is a curved classical law recognized when it is
itself among the candidates?  Both controls retain the separable form
\(x_j(t)=x_{0,j}f(t)\), preserving the ensemble structure of Case A while changing the
temporal~law.

The first control uses the generalized exponential law
\begin{equation}
	x(t)=x_0\exp\!\left(-\mu t^{\beta}\right),
	\label{Eq:genexp}
\end{equation}
with \(\mu=0.15\), \(\beta=0.4\) for stretched relaxation, and~\(\beta=2\) for compressed relaxation; since \(\mu\) multiplies \(t^{\beta}\), it has dimension
\([t]^{-\beta}\), so the same value has a different interpretation for the two values of
\(\beta\).  These controls are false-positive shape tests: in particular, the~stretched
case asks whether a strongly curved but nonfractional decay can appear more
Mittag--Leffler-compatible than a single exponential when the candidate set is
restricted.

The second control uses the rational/Hill-type law
\begin{equation}
	x(t)=\frac{x_0}{1+(t/\mu)^{\beta}},
	\label{Eq:hill}
\end{equation}
with \(\beta=2\) and \(\mu=0.5\).  In~this parameterization, \(\mu\) is a
time-scale parameter rather than the coefficient appearing in the exponent of
Equation~\eqref{Eq:genexp}.  This experiment has a different role from
the generalized-exponential false-positive test: the generating rational/Hill
family is itself included among the fitted candidates.  It is therefore a
matched-family control designed to test whether broadening the candidate set can
recover the appropriate curved classical family, starting from initial guesses that
exclude the generating pair (Section~\ref{Sec:CandidateModels}).

\subsection{Stable kNN Relaxation Profile and Horizon~Extrapolation}
\label{Sec:StableProfile}

All stable benchmarks are evaluated through the same observable, the~kNN relaxation
profile. This subsection defines it and fixes the time coordinate on which the candidate
laws are evaluated.  For~each stable trajectory, the~first \(m_{\mathrm{hist}}=7\)
samples form the history vector, and~future states at horizons \(h=1,\ldots,30\) are the
prediction targets.  Using the trajectory-level split of
Section~\ref{Sec:ScreeningPrinciple}, a~Euclidean \(K\)-nearest-neighbor predictor with
\(K=3\) maps each test history directly to the target at horizon \(h\) by averaging the
corresponding targets of its three nearest training~histories.

The stable observable is a kNN-derived relaxation profile rather than a
forecast-error profile.  Let \(\widehat{x}_i^{(h)}\) denote the direct-horizon
kNN prediction for test trajectory \(i\) at horizon \(h\), and~let
\(x_{i,\mathrm{ref}}\) denote the final state included in its history.  The~profile is
defined as
\begin{equation}\label{Eq:contraction}
	D_{\mathrm{pred}}(h)=
	\frac{\frac{1}{N_{\mathrm{test}}}\sum_i
		\left|\widehat{x}_i^{(h)}\right|}
	{\frac{1}{N_{\mathrm{test}}}\sum_i
		\left|x_{i,\mathrm{ref}}\right|}.
\end{equation}
Thus, \(D_{\mathrm{pred}}(h)\) measures the average amplitude remaining in the
kNN prediction relative to the reference state.  It is a normalized
relaxation-to-equilibrium profile; it is neither a Lyapunov separation between
neighboring trajectories nor a forecast-error~measure.

The time coordinate requires care.  Unlike the exponential, Mittag--Leffler relaxation
does not satisfy a semigroup property, so the clock cannot simply be reset at the end of
the history window while retaining the same scalar relaxation law: the candidate laws must
be evaluated on the absolute target time.
Because the final point of the history window \mbox{occurs at}
\[
t_{\mathrm{ref}}=(m_{\mathrm{hist}}-1)\Delta t=0.6,
\]
the target time is
\begin{equation}\label{Eq:absolute_time}
	t_h=t_{\mathrm{ref}}+\tau_h,
	\qquad
	\tau_h=h\Delta t.
\end{equation}
With \(\Delta t=0.1\), fitting horizons \(h=1,\ldots,20\) therefore correspond to
\(t_h\in[0.7,2.6]\), whereas horizons \(h=21,\ldots,30\) correspond to
\(t_h\in[2.7,3.6]\).

These two intervals are the fitting and horizon-extrapolation intervals of
Section~\ref{Sec:ScreeningPrinciple}.  The~data-separation logic of
Section~\ref{Sec:ScreeningPrinciple} is summarized in
Figure~\ref{Fig:DataSeparationSchematic}.

The normalization in Equation~\eqref{Eq:contraction} affects the profile amplitude but
not the physical time coordinate.  Since every candidate relaxation family contains a free
amplitude parameter \(C\), the~normalization factor is absorbed into \(C\), whereas the
time argument remains the absolute target time \(t_h\).

Because the screening decision is taken on the aggregate profile, we compare the kNN
construction with two simpler constructions, each passed through the same fitting and
extrapolation pipeline. Let \(x_i^{(h)}\) denote the observed state of test trajectory
\(i\) at horizon \(h\), and~let
\(\bar{x}_{\mathrm{ref}}=\frac{1}{N_{\mathrm{test}}}\sum_i|x_{i,\mathrm{ref}}|\) be the
denominator of Equation~\eqref{Eq:contraction}. The~constant training-mean predictor
assigns to every test trajectory the mean training target \(\bar{x}^{(h)}_{\mathrm{train}}\),
giving
\begin{equation}\label{Eq:profile_const}
	D_{\mathrm{const}}(h)=\frac{\bigl|\bar{x}^{(h)}_{\mathrm{train}}\bigr|}{\bar{x}_{\mathrm{ref}}};
\end{equation}
it is an out-of-sample predictor that ignores the individual histories. The~direct
ensemble average,
\begin{equation}\label{Eq:profile_direct}
	D_{\mathrm{direct}}(h)=\frac{1}{\bar{x}_{\mathrm{ref}}}\,
	\frac{1}{N_{\mathrm{test}}}\sum_i\bigl|x_i^{(h)}\bigr|,
\end{equation}
uses the observed futures of the test trajectories and is therefore a descriptive
reference rather than an out-of-sample competitor. To~separate the screening decision from
the accuracy of individual predictions, we also record, for~each of the kNN and constant
predictors, the~normalized individual-trajectory error
\begin{equation}\label{Eq:indiv_error}
	E(h)=\frac{1}{\bar{x}_{\mathrm{ref}}}
	\left[\frac{1}{N_{\mathrm{test}}}\sum_i
	\bigl(x_i^{(h)}-\widehat{x}_i^{(h)}\bigr)^2\right]^{1/2},
\end{equation}
where \(\widehat{x}_i^{(h)}\) denotes the prediction of the predictor being evaluated.
This error measures individual-trajectory prediction accuracy, not the accuracy of the
aggregate screening decision. The~comparison uses the noise-free Cases A--C, the~same
seeds, and~the exponential and Mittag--Leffler~families.

For the observational-noise experiments, the~same profile definition is used with the perturbed observations in both the numerator and the reference-state denominator; the noise protocol is specified in Section~\ref{Sec:NoiseProtocol}.  The~candidate families and their optimization procedures are specified in Section~\ref{Sec:CandidateModels}.

\subsection{Candidate Model Families and~Optimization}
\label{Sec:CandidateModels}

Because a compatibility result is meaningful only relative to the alternatives considered,
each profile is compared with several reference families fitted under a common criterion.
This subsection specifies the objectives, the~candidate families, and~the optimization.
The candidate laws are fitted only to the profile points in the fitting interval.  The~residual objective is defined on the profile itself,
\begin{equation}
	Q_{\mathrm{lin}}
	=
	\frac{1}{n}\sum_{i=1}^{n}
	\left[G_i-\widehat{G}(t_i)\right]^2,
\end{equation}
where \(G_i=D_{\mathrm{pred}}(h_i)\) and \(t_i=t_{h_i}\) is the absolute
target time defined in Equation~\eqref{Eq:absolute_time}.
The stable candidate families are
\begin{align}
	G_{\mathrm{exp}}(t)
	&=C\exp(-\lambda t),\\
	G_{\mathrm{ML}}(t)
	&=C E_{\alpha}\!\left(-\lambda t^{\alpha}\right),\\
	G_{\mathrm{str}}(t)
	&=C\exp(-\lambda t^{\beta}),\\
	G_{\mathrm{RH}}(t)
	&=\frac{C}{1+(t/\mu)^{\beta}},\\
	G_{\mathrm{2exp}}(t)
	&=C_1\exp(-\lambda_1 t)+C_2\exp(-\lambda_2 t),
	\qquad 0<\lambda_1<\lambda_2.
\end{align}
The exponential and Mittag--Leffler models are the common baseline comparison
in Cases A--C.  The~biexponential model is included in the focused noise-free specificity
comparisons of Cases A and B at \(\alpha=0.8\) and Case C at \(s=3\),
as well as in the noisy three-family comparisons of
Section~\ref{Sec:NoiseProtocol}.
The generalized-exponential and rational/Hill families are used in the
separable nonfractional controls described in
Section~\ref{Sec:StableShapeControls}.

All fits were performed using MATLAB R2016a (version 9.0.0.341360) and its \texttt{fminsearch}
implementation of the Nelder--Mead algorithm, which is unconstrained.
Positivity is therefore enforced by exponential parameterizations and finite parameter
intervals by an objective penalty; Table~\ref{tab:candidate_models} lists the
transformations, admissible regions, and~starting values actually used.  In~particular,
the parameterization \(\alpha=\exp(p_3)\) enforces only \(\alpha>0\), and~the restriction
\(0.05\leq\alpha\leq1\) is imposed separately by assigning the objective the value
\(10^{100}\) outside that~interval.

\begin{table}[H]
	\caption{Candidate-model optimization settings.
		All fits use MATLAB \texttt{fminsearch} (Nelder--Mead) with
		\texttt{MaxFunEvals}\(~=5000\), \texttt{MaxIter}\(~=2000\), and~		\texttt{TolFun} = \texttt{TolX}\(~=10^{-10}\). The~solution with the smallest
		finite objective over the listed starts is retained. All fits minimize
		\(Q_{\mathrm{lin}}\).}
	\label{tab:candidate_models}
	\begin{adjustwidth}{0pt}{0pt}
		 
		\footnotesize
		\begin{tabular}{@{}>{\raggedright\arraybackslash}m{2.30cm}>{\raggedright\arraybackslash}m{5.65cm}>{\raggedright\arraybackslash}m{8.65cm}@{}}
			\toprule
			\textbf{Model} & \textbf{Parameterization and Admissible Region} & \textbf{Initialization Actually Used} \\
			\midrule
			Exponential
			& \(C=e^{p_1}\), \(\lambda=e^{p_2}\); \(C,\lambda>0\) by construction.
			& One start:
			\(C_0=\max(G_1,10^{-8})\), \(\lambda_0=0.1\).
			\\ \midrule
			
			Mittag--Leffler
			& \(C=e^{p_1}\), \(\lambda=e^{p_2}\), \(\alpha=e^{p_3}\).
			The log transform enforces only positivity; \(0.05\leq\alpha\leq1\) is imposed
			separately by assigning the objective the value \(10^{100}\) outside the interval.
			& Six \emph{paired} starts, all with
			\(C_0=\max(G_1,10^{-8})\):
			\((\lambda_0,\alpha_0)=
			(0.05,0.2),\allowbreak(0.10,0.5),\allowbreak(0.15,0.7),\allowbreak(0.20,0.85),\allowbreak(0.25,0.95),\allowbreak(0.30,1.0)\).
			\\ \midrule
			
			Generalized exponential
			& \(C=e^{p_1}\), \(\lambda=e^{p_2}\), \(\beta=e^{p_3}\);
			\(0.05<\beta<3\) imposed by objective penalty \(10^{100}\).
			& Three paired starts:
			\((\lambda_0,\beta_0)=(0.05,0.30),(0.10,0.50),(0.20,0.80)\),
			with \(C_0=\max(G_1,10^{-8})\).
			\\ \midrule
			
			Rational/Hill
			& \(C=e^{p_1}\), \(\beta=e^{p_2}\), \(\mu=e^{p_3}\);
			\(0.05<\beta<10\) and \(\mu>0.01\) imposed by objective penalty \(10^{100}\).
			& \(6\times6\) Cartesian grid (36 starts):
			\(\beta_0,\mu_0\in\{0.25,0.4,1.0,1.5,3.0,4.0\}\),
			with \(C_0=\max(G_1,10^{-8})\).
			\\ \midrule
			
			Biexponential
			& \(C_1=e^{p_1}\), \(C_2=e^{p_2}\),
			\(\lambda_1=e^{p_3}\),
			\(\lambda_2=\lambda_1+e^{p_4}\);
			thus \(C_1,C_2>0\) and \(0<\lambda_1<\lambda_2\) hold by construction.
			& \(2\times3\times3\) Cartesian grid (18 starts):
			\(w\in\{0.3,0.7\}\),
			\(\lambda_{1,0}\in\{0.02,0.05,0.1\}\),
			\((\lambda_2-\lambda_1)_0\in\{0.1,0.3,1.0\}\);
			\(C_{1,0}=wG_1\), \(C_{2,0}=(1-w)G_1\).
			\\
			\bottomrule
		\end{tabular}
	\end{adjustwidth}
\end{table}
\vspace{-3pt}

The exponential baseline, a~two-parameter monotone family, uses a single start; the
curved families use multistart fitting, because~different initializations can terminate at
different solutions and objective values.  For~the stable Mittag--Leffler fits in Cases A and B, the~final objectives across the six starts spread by approximately \(4.0\times10^{-4}\) and
\(4.4\times10^{-4}\), while the best solution is many orders of magnitude smaller; this
spread confirms the need for multistart optimization, and~the reported fit is the best
finite-objective solution across starts.  For~the biexponential fits, \(12/18\),
\(9/18\), and~\(3/18\) starts terminated with \texttt{exitflag} = 1 in Cases A, B, and~C,
the others at the iteration/function-evaluation limit; nevertheless, the~final objectives
differed by at most \(7.12\times10^{-8}\), \(7.21\times10^{-8}\), and \mbox{\(7.88\times10^{-10}\).}  Numerically invalid model evaluations receive a large penalty, the~start with the lowest objective is always retained, and~complete per-start diagnostics are
reported in Supplementary Table~S1 and in the archived optimizer logs.  The~Mittag--Leffler function is evaluated with the algorithm of Garrappa
~\cite{Garrappa2015Numerical}.

For model comparison, RMSE is computed on the profile scale.  AIC and BIC use the nominal number of fitted
parameters and serve only as descriptive summaries (Section~\ref{Sec:ScreeningPrinciple});
at parameter boundaries or degenerate limits the effective model dimension may be smaller
than the nominal count.  In~the matched rational/Hill control, where the residuals are at
numerical-precision scale, information-criterion values are not assigned statistical
meaning.

\subsection{Observational-Noise and Uncertainty~Protocol}
\label{Sec:NoiseProtocol}

The benchmarks of Sections~\ref{Sec:StableBenchmarks} and~\ref{Sec:StableProfile} are
noise-free.  To~assess whether the screen remains informative under observational noise,
the analysis is repeated across independent noise realizations rather than relying on a
single noisy run.  This subsection specifies the noise model and the uncertainty
analysis.  Noise is
added to the analytical Case A benchmark (common rate \(\mu=0.15\)) as observation noise
only, without~process noise
in the underlying dynamics: for trajectory \(j\), the~clean samples are perturbed
according to
\begin{equation}
	x^{\mathrm{obs}}_j(t_n)
	=
	x_j(t_n)+\varepsilon_{j,n},
	\qquad
	\varepsilon_{j,n}\overset{\mathrm{i.i.d.}}{\sim}
	\mathcal{N}(0,\sigma_j^2),
	\qquad
	\sigma_j=0.03\,|x_{0,j}|.
	\label{Eq:noise_model}
\end{equation}
Thus, the~noise standard deviation is fixed in time along each trajectory and
equals \(3\%\) of that trajectory's initial amplitude; it is not \(3\%\) of the
instantaneous state magnitude.  Noise samples are independent across time and
across~trajectories.

For each random seed, the~trajectory-level training/test partition is generated
before the noise realization is added, after~which the complete ensemble is
perturbed.  Consequently, both the kNN training data and the kNN test data
contain noisy observations.  The~relaxation profile retains the definition in
Equation~\eqref{Eq:contraction}, but~in the noisy experiment both the predicted
future amplitudes in the numerator and the reference-state amplitudes in the
denominator are formed from the perturbed observations.  The~same
training-only standardization and direct-horizon kNN construction are then
applied as in the noise-free stable~analysis.

The uncertainty study uses \(50\) independent random seeds for each generating
order \(\alpha\in\{0.6,0.7,0.8,0.9,1\}\).  The~runs are paired across
\(\alpha\): for a given seed, the~initial conditions, trajectory partition, and~additive noise realization are kept the same, while only the clean generating
trajectory changes with \(\alpha\).  Consequently, differences across generating
orders cannot be attributed to different ensembles, train/test partitions, or~noise realizations; the pairing isolates the change associated with the clean
relaxation law as \(\alpha\) varies.

For each seed we record the fitted order \(\alpha_{\mathrm{fit}}\) and the
horizon-extrapolation difference
\begin{equation}
	\Delta\mathrm{RMSE}
	=
	\mathrm{RMSE}^{*}_{\mathrm{exp}}
	-
	\mathrm{RMSE}^{*}_{\mathrm{ML}},
	\label{Eq:delta_rmse}
\end{equation}
where the superscript \(*\) denotes evaluation on horizons \(21\)--\(30\)
without refitting.  Positive values indicate a smaller
horizon-extrapolation RMSE for the Mittag--Leffler model.  Across the \(50\)
paired runs, we report the mean together with the \(95\%\) confidence
interval for the mean, based on the Student \(t\) distribution with 49 degrees
of freedom, of~both \(\Delta\mathrm{RMSE}\) and
\(\alpha_{\mathrm{fit}}\), and~separately count wins, losses, and~exact ties.  Operationally, a~tie is
recorded when \(|\Delta\mathrm{RMSE}|<10^{-8}\) or when \(\alpha_{\mathrm{fit}}\) lies
within \(10^{-6}\) of 1, where the Mittag--Leffler model reduces numerically to the
exponential~model.

To complement the across-realization analysis with an estimate of test-ensemble
variability conditional on a fixed fitted predictor, a~trajectory bootstrap is performed
for one representative noisy run.  The~fitted kNN training set is held fixed, while the test trajectories are
resampled with replacement.  For~each of \(B=500\) bootstrap resamples, the~noisy
relaxation profile is recomputed, the~exponential and Mittag--Leffler models are
refitted on horizons 1--20, and~\(\Delta\mathrm{RMSE}\) is reevaluated on
horizons 21--30.  The~reported bootstrap interval is the percentile \(95\%\)
interval of these \(500\) values.  It therefore quantifies uncertainty
conditional on the estimated training predictor and the selected representative
run.  It is not interpreted as a bootstrap confidence interval for the full
pipeline, because~the training trajectories, trajectory split, and~model-building
stage are not~resampled.

The noise-free specificity analysis showed that a classical distributed-rate ensemble can
favor the Mittag--Leffler law over a single exponential, and~that this preference
disappears once a biexponential alternative is admitted
(Section~\ref{Sec:SpecificityResults}). To~test whether this conclusion survives
observational noise, the~same noise model, the~same \(50\) seeds, and~the same fitting and
extrapolation horizons are applied to two classical generators: Case B with \(\alpha=1\),
whose trajectories are exponentials with rates \(\mu_j\sim\mathcal{U}[0.05,0.25]\), and~Case C with \(s=3\). Each noisy profile is compared first with the exponential and
Mittag--Leffler families, as~above, and~then with the exponential, Mittag--Leffler, and~biexponential families. A~broader candidate set that rejects the Mittag--Leffler law for
classical data could also do so for fractional data. To~check that the three-family
comparison retains sensitivity to fractional relaxation, it is also applied to the
fractional Case B with \(\alpha=0.8\); this single reference is a targeted check, not a
general assessment of~sensitivity.

The selection rule was fixed before the results were examined. For~each seed, the~family
with the lowest horizon-extrapolation RMSE is selected, and~families whose RMSE differs
from the lowest by less than \(10^{-8}\) are counted as tied. The~convention that a
Mittag--Leffler fit with \(|\alpha_{\mathrm{fit}}-1|<10^{-6}\) coincides with the
exponential applies only to that pair and does not create a tie with the biexponential.
Selection frequencies are reported with \(95\%\) Wilson score intervals, and~the paired
held-out RMSE differences between families with \(95\%\) Student \(t\) confidence intervals
for their means. Seeds with \(|\alpha_{\mathrm{fit}}-1|<10^{-6}\) are also counted as
boundary hits, separately from ties defined by RMSE. For~the binary comparisons, exact
sign tests and Wilcoxon signed-rank tests on the paired differences, excluding ties under
the stated binary rule, are reported as supporting, unadjusted analyses. Because~the
biexponential is the most flexible family, its fits are monitored for convergence and
degeneracy (a normalized component weight or relative rate gap below \(10^{-3}\)); these
diagnostics are reported in the Supplementary~Material.

\subsection{Stable Sensitivity~Protocol}
\label{Sec:SensitivityProtocol}

The practical value of the screen also depends on its robustness to reasonable analysis
choices.  To~assess this, the~stable pipeline is rerun on Case B at the representative generating
order \(\alpha=0.8\), changing one analysis setting at a time relative to the baseline
configuration
\[
K=3,\qquad
m_{\mathrm{hist}}=7,\qquad
\Delta t=0.1,\qquad
\text{fitting horizons }1\text{--}20.
\]
The four one-at-a-time perturbations are
\begin{align}
	&K\in\{1,3,5,10\},\\
	&m_{\mathrm{hist}}\in\{4,7,10\},\\
	&\text{fitting horizons}\in\{1\text{--}15,\;1\text{--}20,\;1\text{--}25\},\\
	&\Delta t\in\{0.05,0.1,0.2\}.
\end{align}
This design provides a local sensitivity assessment around the baseline configuration
rather than an exhaustive exploration of interactions among analysis choices.
Case B is used because, unlike Case A, it is not reduced to a single temporal template by
separability.  The~same underlying data-generation setting and the same ten random seeds
are used across configurations, so each row is compared against the same baseline
experiment.

For the sampling-interval perturbation, the~physical history, forecast, and~fitting windows are kept fixed rather than keeping the number of samples fixed.  Consequently, changing \(\Delta t\) also changes the discrete
history length and the number of forecast/fitting horizons.  Specifically,
\[
\Delta t=0.05:
\quad
m_{\mathrm{hist}}=13,\quad
H_{\max}=60,\quad
\text{fitting horizons }1\text{--}40,
\]
whereas
\[
\Delta t=0.2:
\quad
m_{\mathrm{hist}}=4,\quad
H_{\max}=15,\quad
\text{fitting horizons }1\text{--}10.
\]
These choices preserve the baseline physical windows exactly: the history span
is \mbox{\((m_{\mathrm{hist}}-1)\Delta t=0.6\),} the~maximum forecast span is
\(H_{\max}\Delta t=3.0\), and~the fitting span is \(2.0\) time units in all
three sampling configurations.  Thus, the~sampling-interval rows are not
independent perturbations of \(m_{\mathrm{hist}}\): they probe a change in
temporal resolution while preserving the physical~windows.

For each configuration, the~same stable fitting procedure described in
Sections~\ref{Sec:StableProfile} and~\ref{Sec:CandidateModels} is applied, and~the fitted order \(\alpha_{\mathrm{fit}}\), the~horizon-extrapolation
\(\Delta\mathrm{RMSE}\) of \mbox{Equation~\eqref{Eq:delta_rmse}}, and~the descriptive
information-criterion difference are summarized over the repeated seeds.  Since
the number of fitted profile points changes when \(\Delta t\) or the fitting
range changes, AIC differences are interpreted only within configurations with
the same number of fitted points and are not compared numerically across all
sensitivity~rows.

\subsection{Numerical Reliability of a Single-Record~Example}
\label{Sec:NumericalCheck}

A screen applied to simulated observations can be no more reliable than the simulation that
produced them. This matters for the single-record example, a~scalar series generated from a
nonlinear fractional system initially investigated as a chaotic test case. The~system is
the fractional system analyzed in~\cite{sene2021analysis},
\begin{equation}
	\begin{cases}
		{}_{C}D_{0,t}^{\alpha}x=ax-by-yz,\\
		{}_{C}D_{0,t}^{\alpha}y=cx,\\
		{}_{C}D_{0,t}^{\alpha}z=-dz+y^{2},
	\end{cases}
	\label{Eq:sene_system}
\end{equation}
with \(a=-2\), \(b=-6.4\), \(c=1\), and~\(d=1\); related Caputo fractional-order chaotic
systems have also been investigated in~\cite{Sene2022Theory}. It is integrated to \(T=200\)
from \(\mathbf{x}_0=(0.2,0.2,0.2)^{\mathsf T}\) with \(\alpha=0.916\), and~the
\(x\)-component is observed after the first \(40\) time units are discarded as~transient.

The record was originally generated with the explicit product-integration (fractional
Euler) scheme~\cite{Diethelm2005Algorithms} at step \(\Delta t=0.01\). To~verify whether the
regime it produces is robust, we use Garrappa's MATLAB implementation \texttt{fde12} of the
Adams--Bashforth--Moulton method~\cite{Diethelm2005Algorithms}, both in predictor-only mode,
which reproduces the original explicit discretization, and~in predictor--corrector (PECE)
mode, with~\(\Delta t\in\{0.01,0.005,0.0025,0.00125\}\). To~assess whether the observed
regime persists under numerical refinement, we compare both summary diagnostics and
differences between numerical solutions: the standard deviation of \(x\) over \(t\geq40\)
and over the final window \(180\leq t\leq200\), the~distance of the final state from the
nearest nonzero equilibrium, and~the root-mean-square differences between PECE solutions at
successive step sizes on a common grid. A~run with an iterated corrector at
\(\Delta t=0.005\) and a control run at \(\alpha=0.93\) complete the~check.

These diagnostics are complemented by the local stability of the nonzero equilibria
\(E_\pm=(0,\pm\sqrt{6.4},6.4)\). For~a commensurate Caputo system, an~equilibrium is locally
asymptotically stable if every eigenvalue \(\lambda_j\) of the Jacobian satisfies
\(|\arg\lambda_j|>\alpha\pi/2\) \cite{diethelm2010analysis}. This criterion characterizes
local behavior only and does not determine the global dynamics from a given initial
condition. The~kNN screening protocol originally applied to this record is documented in
the Supplementary~Material.

\section{Results}
\label{Sec:Results}

The results progress from matched-family calibration through increasingly demanding
tests of heterogeneity and specificity, followed by comparisons of profile
constructions, observational noise, and~analysis sensitivity. The~section concludes
with a numerical-reliability check of the single-record~example.

\subsection{Matched-Family~Calibration}
\label{Sec:ExactCalibration}

\begin{figure}[H]
	\begin{adjustwidth}{0pt}{0pt}
		\centering
		\setlength{\tabcolsep}{3pt}%
		\begin{tabular}{@{}ccc@{}}
			\includegraphics[width=0.32\linewidth]{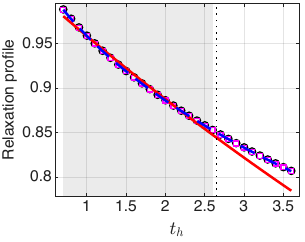} &
			\includegraphics[width=0.32\linewidth]{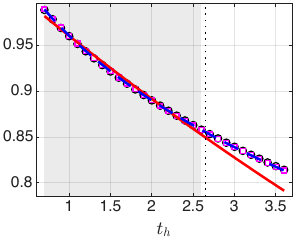} &
			\includegraphics[width=0.32\linewidth]{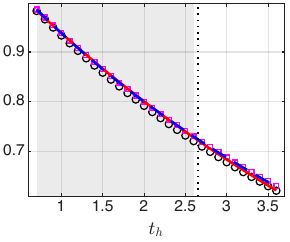} \\[2pt]
			
			\includegraphics[width=0.32\linewidth]{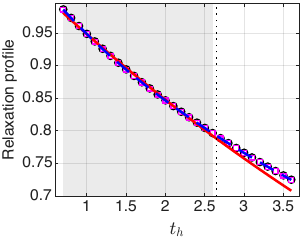} &
			\includegraphics[width=0.32\linewidth]{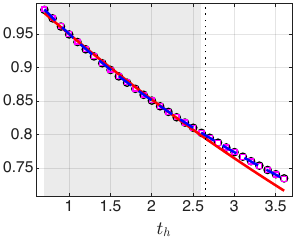} &
			\includegraphics[width=0.32\linewidth]{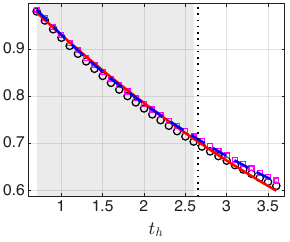} \\[2pt]
			
			\includegraphics[width=0.32\linewidth]{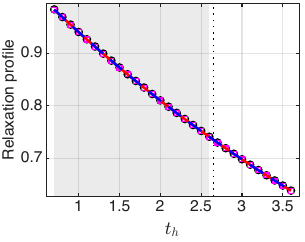} &
			\includegraphics[width=0.32\linewidth]{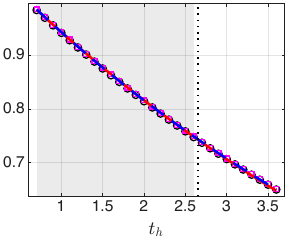} &
			\includegraphics[width=0.32\linewidth]{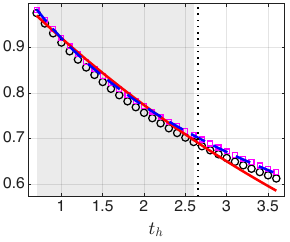} \\[2pt]
			
			(\textbf{a}) \textbf{Case A} &
			(\textbf{b}) \textbf{Case B} &
			(\textbf{c}) \textbf{Case C} \\[4pt]
			
			\multicolumn{3}{c}{%
			\includegraphics[width=0.95\linewidth]{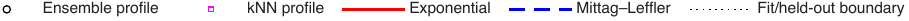}}
		\end{tabular}
	\end{adjustwidth}
	
	\caption{Stable relaxation profiles for the three benchmark cases. Each
		panel shows the directly computed ensemble profile and the kNN-predicted
		profile, together with the exponential and Mittag--Leffler laws fitted on
		horizons $h=1$--$20$ (shaded) and continued over the held-out horizons
		$h=21$--$30$, on~the absolute target time
		$t_h=t_{\mathrm{ref}}+h\,\Delta t$.
		(\textbf{a})~Case A, homogeneous rate $\mu=0.15$;
		(\textbf{b})~Case B, heterogeneous rates $\mu_j\sim\mathcal{U}[0.05,0.25]$;
		in both columns, from~top to bottom, $\alpha=0.6$, $0.8$ and $1$.
		(\textbf{c})~Case C, purely exponential trajectories with log-uniform rates
		(Equation~\eqref{Eq:case_c_rate}); from top to bottom, $s=3$, $5$ and $10$.}
	\label{Fig:StableCore}
\end{figure}

Case A provides a matched, noise-free calibration of the complete stable
pipeline.  Because~the trajectories are generated analytically in an
exactly matched setting, the~fitted-order discrepancy is observed at the
\(10^{-11}\) scale.  Across
\(\alpha\in\{0.6,0.7,0.8,0.9\}\), the~fitted Mittag--Leffler order therefore
reproduces the generating value to numerical precision, and~the fitted rate is
\(\lambda_{\mathrm{ML}}=0.150000\) in every case.  The~corresponding
Mittag--Leffler horizon-extrapolation residuals are likewise at
numerical-precision scale (Table~\ref{tab:exact_calibration}); they are not
interpreted as evidence of statistical discrimination.  At~\(\alpha=1\), the~fitted order reaches the exponential limit and the two families coincide to
numerical precision. Representative profiles for Cases A--C are shown in
Figure~\ref{Fig:StableCore}.

\begin{table}[H]
	\caption{Exact matched calibration (Case A).  Values of
		\(\alpha_{\mathrm{fit}}\), RMSE, \(\Delta\mathrm{AIC}\), and~		\(\Delta\mathrm{BIC}\) are shown for the first seed; the maximum absolute order
		deviation is taken over ten independent seeds.  The~seed-to-seed dispersion of
		\(\alpha_{\mathrm{fit}}\) is at the \(10^{-11}\) scale.
		RMSE\(^{*}\) is evaluated on horizons \(21\)--\(30\) without refitting.
		Negative \(\Delta\mathrm{AIC}\) and \(\Delta\mathrm{BIC}\) favor the
		Mittag--Leffler model.}
	\label{tab:exact_calibration}

\begin{adjustwidth}{0pt}{0pt}
\begin{tabularx}{\fulllength}{CCCCCCC}
			\toprule
			\boldmath{\(\alpha\) }&
			\boldmath{\(\alpha_{\mathrm{fit}}\) }&
			\boldmath{\(\max|\alpha_{\mathrm{fit}}-\alpha|\) }&
			\boldmath{\textbf{RMSE}\(^{*}_{\mathrm{exp}}\) }&
			\boldmath{\textbf{RMSE}\(^{*}_{\mathrm{ML}}\) }&
			\boldmath{\(\Delta\mathrm{AIC}\) }&
			\boldmath{\(\Delta\mathrm{BIC}\) }\\
			\midrule
			0.6 & 0.600000 & \(1.5\times10^{-11}\) &
			\(1.485\times10^{-2}\) & \(5.393\times10^{-14}\) &
			\(-964.07\) & \(-963.07\) \\
			0.7 & 0.700000 & \(1.0\times10^{-11}\) &
			\(1.388\times10^{-2}\) & \(2.925\times10^{-13}\) &
			\(-919.13\) & \(-918.14\) \\
			0.8 & 0.800000 & \(2.2\times10^{-11}\) &
			\(1.124\times10^{-2}\) & \(2.040\times10^{-13}\) &
			\(-922.00\) & \(-921.01\) \\
			0.9 & 0.900000 & \(2.7\times10^{-11}\) &
			\(6.683\times10^{-3}\) & \(1.316\times10^{-12}\) &
			\(-889.10\) & \(-888.10\) \\
			1.0 & 1.000000 & \(1.2\times10^{-10}\) &
			\(2.694\times10^{-12}\) & \(5.295\times10^{-13}\) &
			-- & -- \\
			\bottomrule
		\end{tabularx}
\end{adjustwidth}

	\noindent{\footnotesize
		At \(\alpha=1\), the~Mittag--Leffler family reduces to the exponential
		special case.  Both residuals are at numerical-precision scale, so
		information-criterion differences are not interpreted.}
\end{table}
\vspace{-4pt}

Because the generator and fitted family coincide and the Case A ensemble is
multiplicatively separable, this near-exact recovery is treated as a
matched-family calibration result rather than as evidence of general
fractional-order~identifiability.

\subsection{Heterogeneous Relaxation~Rates}
\label{Sec:HeterogeneousRates}

Case B breaks the multiplicative separability of Case A by assigning each
trajectory its own relaxation rate
\(\mu_j\sim\mathcal{U}[0.05,0.25]\) while keeping a common generating order.
The fitted Mittag--Leffler order remains close to the generating value across
the sweep, but~with a small systematic negative offset of approximately
\(0.02\) (Table~\ref{tab:heterogeneous_rates}).  The~fitted rate remains close
to the mean generating rate, with~\(\lambda_{\mathrm{ML}}\approx0.149\).  Once separability is removed, the~fitted order
in this heterogeneous ensemble thus behaves as an effective profile parameter rather than reproducing the generating
order~exactly.

The Mittag--Leffler model has lower horizon-extrapolation error than the
single exponential for every generating order in this heterogeneous ensemble.  Notably, this remains true at \(\alpha=1\), where every individual
trajectory is exponential but the ensemble contains a distribution of rates:
the fitted order is \(0.9768\pm0.0009\) and
\(\Delta\mathrm{AIC}=-94.6\pm0.4\).  This integer-order case anticipates the
more stringent classical distributed-rate specificity test in
Section~\ref{Sec:SpecificityResults}.

 \begin{table}[H]
	\caption{Case B with heterogeneous relaxation rates
		\(\mu_j\sim\mathcal{U}[0.05,0.25]\).  RMSE values are for the first random
		seed; \(\alpha_{\mathrm{fit}}\) and \(\Delta\mathrm{AIC}\) are mean
		\(\pm\) standard deviation over ten seeds.
		RMSE\(^{*}\) is evaluated on horizons \(21\)--\(30\) without refitting.
		Negative \(\Delta\mathrm{AIC}\) favors the Mittag--Leffler model.}
	\label{tab:heterogeneous_rates}
	\begin{adjustwidth}{0pt}{0pt}
	
		\begin{tabularx}{\fulllength}{CCCCCC}
			\toprule
			\boldmath{\(\alpha\)} &
			\boldmath{\(\alpha_{\mathrm{fit}}\)} &
			\boldmath{\(\alpha_{\mathrm{fit}}-\alpha\)} &
			\textbf{RMSE}\boldmath{\(^{*}_{\mathrm{exp}}\)} &
			\textbf{RMSE}\boldmath{\(^{*}_{\mathrm{ML}}\) }&
			\boldmath{\(\Delta\mathrm{AIC}\)} \\
			\midrule
			0.6 & \(0.5828\pm0.0010\) & \(-0.017\) &
			\(1.496\times10^{-2}\) & \(5.508\times10^{-5}\) &
			\(-245.2\pm2.8\) \\
			0.7 & \(0.6805\pm0.0010\) & \(-0.020\) &
			\(1.423\times10^{-2}\) & \(9.654\times10^{-5}\) &
			\(-222.4\pm2.5\) \\
			0.8 & \(0.7789\pm0.0010\) & \(-0.021\) &
			\(1.202\times10^{-2}\) & \(1.606\times10^{-4}\) &
			\(-198.3\pm2.0\) \\
			0.9 & \(0.8777\pm0.0010\) & \(-0.022\) &
			\(7.933\times10^{-3}\) & \(2.390\times10^{-4}\) &
			\(-167.0\pm1.7\) \\
			1.0 & \(0.9768\pm0.0009\) & \(-0.023\) &
			\(1.882\times10^{-3}\) & \(3.449\times10^{-4}\) &
			\(-94.6\pm0.4\) \\
			\bottomrule
		\end{tabularx}
	\end{adjustwidth}
\end{table}
\vspace{-12pt}

\subsection{Specificity Tests and Candidate-Set Comparisons}
\label{Sec:SpecificityResults}

The classical distributed-rate control in Case C first exposes a limitation of
the binary exponential-versus--Mittag--Leffler comparison.  All trajectories are
strictly exponential, but~their rates are drawn from the log-uniform family
defined in Section~\ref{Sec:StableABC}.  Nevertheless, the~Mittag--Leffler model
is preferred over the single exponential in all ten runs at each tested spread
factor (Table~\ref{tab:distributed_rate_binary}).  As~the rate distribution
widens from \(s=3\) to \(s=10\), the~fitted order decreases monotonically from
\(0.933\) to \(0.668\).  The~separation between the exponential and
Mittag--Leffler continuations over the held-out horizons grows with \(s\)
(Figure~\ref{Fig:StableCore}c).

\begin{table}[H]
	\caption{Case C: binary exponential-versus--Mittag--Leffler comparison for
		classical distributed-rate ensembles.  The~trajectories are purely
		exponential, with~\(\mu_j=0.15s^{u_j}\),
		\(u_j\sim\mathcal{U}[-1,1]\).  \(\Delta\mathrm{RMSE}=
		\mathrm{RMSE}^{*}_{\mathrm{exp}}-\mathrm{RMSE}^{*}_{\mathrm{ML}}\), in~units
		of \(10^{-3}\); positive values indicate lower horizon-extrapolation RMSE for
		the Mittag--Leffler fit.  Values are mean \(\pm\) standard deviation over ten
		seeds.}
	\label{tab:distributed_rate_binary}

		\begin{tabularx}{\textwidth}{CCCCC}
			\toprule
			\boldmath{\(s\)} & \textbf{Rate Range} & \boldmath{\(\alpha_{\mathrm{fit}}\)} &
			\boldmath{\(\Delta\mathrm{RMSE}\)} & \textbf{ML Wins} \\
			\midrule
			3  & \([0.05,0.45]\)   & \(0.933\pm0.003\) &
			\(5.36\pm0.35\)  & 10/10 \\
			5  & \([0.03,0.75]\)   & \(0.847\pm0.008\) &
			\(12.53\pm1.27\) & 10/10 \\
			10 & \([0.015,1.5]\)   & \(0.668\pm0.022\) &
			\(26.43\pm3.20\) & 10/10 \\
			\bottomrule
		\end{tabularx}
	 
\end{table}

The same comparison changes when a more flexible classical relaxation family is
included.  Table~\ref{tab:expanded_candidate_comparison} compares the
exponential, Mittag--Leffler, and~biexponential fits for Case A and Case B at
\(\alpha=0.8\), and~for Case C at \(s=3\).  In~Case B, the~Mittag--Leffler model remains the lowest-error model on
the horizon-extrapolation interval in all ten runs.  Case A shows the same
ranking, but~it is the matched-family, noise-free calibration of
Section~\ref{Sec:ExactCalibration}, so its machine-precision Mittag--Leffler
residuals are not interpreted as a statistically meaningful advantage over the
biexponential model.  In~Case C, by~contrast, the~biexponential model is best
in all ten runs: its mean held-out RMSE is \(0.0129\times10^{-3}\), compared
with \(0.992\times10^{-3}\) for the Mittag--Leffler model and
\(6.36\times10^{-3}\) for the single exponential.  Thus, in~the noise-free setting, the~preference
observed in the binary Case C comparison does not persist after expansion of
the candidate set.  For~Case C, the~apparent Mittag--Leffler preference is therefore a property
of this restricted candidate comparison, not evidence that the classical ensemble is fractionally
generated.  Section~\ref{Sec:NoiseResults} examines whether this conclusion survives
observational~noise.

The separable nonfractional shape controls give two complementary checks.  For~the stretched-exponential generator (\(\beta=0.4\)), the~Mittag--Leffler fit
produces an interior fitted order below one and markedly improves on the single
exponential.  For~the compressed control (\(\beta=2\)), the~fitted order reaches
\(\alpha_{\mathrm{fit}}=1\), so the Mittag--Leffler model reduces to the
exponential limit; however, neither family describes this control well, with~training RMSE \(2.925\times10^{-2}\) and horizon-extrapolation RMSE
\(9.913\times10^{-2}\) for both.  The~boundary solution therefore reflects
failure of the tested exponential/Mittag--Leffler families to capture this
compressed shape rather than successful~discrimination.

\begin{table}[H]
	\caption{Expanded candidate comparison for Case A and Case B at
		\(\alpha=0.8\), and~Case C at \(s=3\).  Held-out RMSE values are in units of
		\(10^{-3}\) and are mean \(\pm\) standard deviation over ten seeds.  ``Wins''
		gives the number of runs in which the model has the lowest held-out RMSE.
		Case A is included only as a matched-family calibration reference; its
		machine-precision Mittag--Leffler residuals are interpreted as in
	 Section~\ref{Sec:ExactCalibration}, not as evidence of statistical
		discrimination.}
	\label{tab:expanded_candidate_comparison}
	 \footnotesize
	 
		\begin{tabularx}{\textwidth}{CLCC}
			\toprule
			\textbf{Case} & \textbf{Model} & \textbf{Held-Out RMSE} & \textbf{Wins} \\
			\midrule
			\multirow{3}{*}{A} & Exponential       & \(11.242\pm0.002\) & 0/10 \\
			& Mittag--Leffler   & \(<10^{-9}\)     & 10/10 \\
			& Biexponential     & \(0.7285\pm0.0002\)& 0/10 \\
			\midrule
			\multirow{3}{*}{B} & Exponential       & \(12.07\pm0.32\) & 0/10 \\
			& Mittag--Leffler   & \(0.181\pm0.011\)& 10/10 \\
			& Biexponential     & \(0.764\pm0.022\)& 0/10 \\
			\midrule
			\multirow{3}{*}{C} & Exponential       & \(6.36\pm0.43\)  & 0/10 \\
			& Mittag--Leffler   & \(0.992\pm0.079\)& 0/10 \\
			& Biexponential     & \(0.0129\pm0.0013\)& 10/10 \\
			\bottomrule
		\end{tabularx}
 
\end{table}

The rational/Hill control serves a different purpose because the generating
family itself is included among the candidates.  The~fit recovers
\(\beta_{\mathrm{fit}}=2.000000\) and
\(\mu_{\mathrm{fit}}=0.500000\), with~training and horizon-extrapolation RMSEs
of \(4.764\times10^{-13}\) and \(1.279\times10^{-13}\), respectively.  These
machine-precision residuals are treated only as a same-family numerical
recovery check; information-criterion differences at this scale are not
assigned statistical~meaning.

\subsection{Contribution of the kNN~Construction}
\label{Sec:BaselineResults}

Because the screening decision is taken on the aggregate profile, the~kNN profile was
compared with the constant training-mean predictor and the direct ensemble average
(Section~\ref{Sec:StableProfile}) in all thirteen noise-free configurations of Cases A--C,
with ten seeds each. The~three constructions lead to the same
exponential-versus--Mittag--Leffler outcome in every seed of every configuration
(Table~\ref{tab:baselines}). For~example, in~Case A at \(\alpha=1\), all ten seeds yield
ties under each construction. Their fitted orders are also close: the largest difference in
the mean fitted order, \(0.015\), occurs in Case C at \(s=5\), and~in no configuration does
it change the screening outcome. The~constructions differ sharply, however, at~the level of
individual trajectories. The~mean held-out individual error of the kNN predictor ranges from
\(0.002\) in Case A to \(0.077\) in Case C at \(s=10\), whereas that of the constant
predictor ranges from \(0.29\) to \(0.50\), between~about \(5.6\) and \(135\) times larger.
Part of this advantage reflects the amplitude information contained in each observed
history, which the constant predictor ignores. The~complete summary for all thirteen
configurations is given in the Supplementary~Material.

\begin{table}[H]
	\caption{\textls[-15]{Comparison of profile constructions in representative noise-free
		configurations (ten seeds each). \(\alpha_{\mathrm{fit}}\): mean fitted order for the
		kNN, constant training-mean (Const), and~direct ensemble-average (Direct) profiles.
		W/L/T: Mittag--Leffler wins, losses, and~ties on the held-out horizons, identical for
		the three constructions. \(E\): mean held-out individual-trajectory error,
		Equation~\eqref{Eq:indiv_error}, for~the kNN and constant predictors; it is not
		defined for the direct average, which uses the observed futures.}}
	\label{tab:baselines}
	 
	\begin{tabularx}{\textwidth}{LlCCCCc}
		\toprule
		& & \multicolumn{3}{c}{\boldmath{\(\alpha_{\mathrm{fit}}\)}} & & \boldmath{\(E\)} \\
		\cmidrule{3-5}  \cmidrule{7-7}
		\textbf{\multirow{-2}{*}{Case}} & \textbf{\multirow{-2}{*}{Parameter}} & \textbf{kNN} & \textbf{Const} & \textbf{Direct} & \textbf{\multirow{-2}{*}{W/L/T}} & \textbf{kNN/Const} \\
		\midrule
		A & \(\alpha=0.8\) & 0.800 & 0.800 & 0.800 & 10/0/0 & 0.0024 / 0.318 \\
		B & \(\alpha=0.8\) & 0.779 & 0.775 & 0.775 & 10/0/0 & 0.040 / 0.339 \\
		B & \(\alpha=1\)   & 0.977 & 0.972 & 0.973 & 10/0/0 & 0.053 / 0.320 \\
		C & \(s=3\)        & 0.933 & 0.924 & 0.926 & 10/0/0 & 0.063 / 0.354 \\
		C & \(s=5\)        & 0.847 & 0.832 & 0.836 & 10/0/0 & 0.070 / 0.416 \\
		C & \(s=10\)       & 0.668 & 0.659 & 0.667 & 10/0/0 & 0.077 / 0.499 \\
		\bottomrule
	\end{tabularx}
\end{table}
\unskip

\subsection{Observational Noise and~Uncertainty}
\label{Sec:NoiseResults}

The next question is whether the matched-family screening behavior of Case A persists
once observational noise is introduced.  With~trajectory-scaled noise, the~fitted order remains close to the
generating value over the fractional cases, while the horizon-extrapolation
advantage of the Mittag--Leffler fit decreases as \(\alpha\) approaches the
integer-order limit (Table~\ref{tab:noise_results}).  In~contrast to Case B,
no comparable systematic negative offset is observed here across
\(\alpha=0.6\)--\(0.9\): the corresponding estimates are
\(0.607\), \(0.706\), \(0.805\), and~\(0.904\).  For~\(\alpha=0.6,0.7,0.8,\) and \(0.9\), the~95\% confidence interval for
\(\Delta\mathrm{RMSE}=
\mathrm{RMSE}^{*}_{\mathrm{exp}}-\mathrm{RMSE}^{*}_{\mathrm{ML}}\)
remains above zero across the fifty-seed experiment.  At~\(\alpha=1\), the~mean difference is slightly negative,
\(-0.39\times10^{-3}\), with~a 95\% confidence interval
\([-0.74,-0.03]\times10^{-3}\); 30 of the 50 runs are exact numerical ties,
15 favor the exponential model, and~5 favor the Mittag--Leffler model.  All 30 ties at
\(\alpha=1\) are boundary hits, with~\(\alpha_{\mathrm{fit}}\) within \(10^{-6}\) of 1; at
\(\alpha=0.9\) a single seed reaches the boundary, and~none does at lower~orders.

\begin{table}[H]
	\caption{Noisy stable experiment with
		\(\sigma_j=0.03|x_{0,j}|\).  Values are based on fifty independent seeds.
		\(\Delta\mathrm{RMSE}=
		\mathrm{RMSE}^{*}_{\mathrm{exp}}-\mathrm{RMSE}^{*}_{\mathrm{ML}}\) is evaluated
		on horizons \(21\)--\(30\) without refitting and is reported in units of
		\(10^{-3}\).  Positive values indicate lower horizon-extrapolation RMSE for the
		Mittag--Leffler model. Bracketed values are 95\% confidence intervals for the
		mean across the fifty realizations.}
	\label{tab:noise_results}

		\begin{tabularx}{\textwidth}{CCC}
			\toprule
			\boldmath{\(\alpha\)} &
			\boldmath{\(\alpha_{\mathrm{fit}}\)} \textbf{[95\% CI]} &
			\boldmath{\(\Delta\mathrm{RMSE}\)} \textbf{[95\% CI]} \\
			\midrule
			0.6 & \(0.607\;[0.585,0.629]\) & \(11.48\;[10.82,12.15]\) \\
			0.7 & \(0.706\;[0.689,0.723]\) & \(10.49\;[9.81,11.16]\) \\
			0.8 & \(0.805\;[0.791,0.820]\) & \(7.87\;[7.22,8.53]\) \\
			0.9 & \(0.904\;[0.892,0.916]\) & \(3.48\;[2.82,4.14]\) \\
			1.0 & \(0.987\;[0.981,0.993]~^{\dagger}\) &
			\(-0.39\;[-0.74,-0.03]~^{\ddagger}\) \\
			\bottomrule
		\end{tabularx}

	\noindent{\footnotesize
		\({}^{\dagger}\)~The fitted order is constrained by
		\(\alpha_{\mathrm{fit}}\leq 1\), and~30 of the 50 fits at \(\alpha=1\) lie on this
		boundary. The~reported Student \(t\) interval summarizes the across-seed mean of the
		constrained fitted orders; it should not be interpreted as a confidence interval for
		the generating order or for an unconstrained estimator.
		\({}^{\ddagger}\)~For \(\alpha=1\), the~50 runs comprise 5 Mittag--Leffler
		wins, 15 exponential wins, and~30 exact numerical ties
		(\(|\Delta\mathrm{RMSE}|<10^{-8}\) or \(|\alpha_{\mathrm{fit}}-1|<10^{-6}\)).}
\end{table}

\begin{figure}[H]
	 
	\hspace{-9pt}\includegraphics[width=0.80\linewidth]{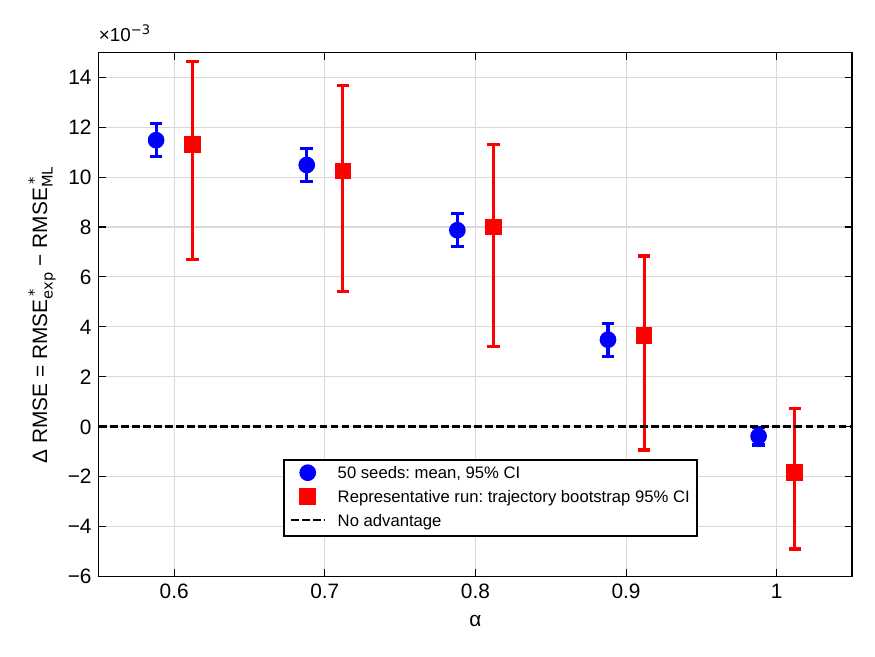}
	\caption{Held-out difference \(\Delta\mathrm{RMSE}=\mathrm{RMSE}^{*}_{\mathrm{exp}}
		-\mathrm{RMSE}^{*}_{\mathrm{ML}}\) versus the generating order for the noisy
		stable experiment.  The~first series shows the mean and the 95\% confidence interval for the mean
		across fifty independent realizations.  The~second shows the
		percentile 95\% trajectory-bootstrap interval for one representative
		realization, conditional on its estimated kNN predictor and fixed training
		set.  The~two intervals therefore quantify different sources of uncertainty
		and should not be read as competing estimates of the same sampling
		distribution.}
	\label{Fig:NoiseUncertainty}
\end{figure}

Figure~\ref{Fig:NoiseUncertainty} displays the same held-out comparison
together with the trajectory-bootstrap interval for one representative run.  The~two
intervals quantify different sources of uncertainty: the first summarizes
variation across fifty independently generated realizations, whereas the
bootstrap interval is conditional on the estimated kNN predictor and fixed
training set of a single representative realization.  In~the latter calculation,
test trajectories are resampled with replacement and the profile and model fits
are recomputed for each bootstrap~sample.

The noise-free specificity result of Section~\ref{Sec:SpecificityResults} was then tested
under the same noise model. In~the binary comparison, the~classical controls frequently
favor the Mittag--Leffler law. In~Case C with \(s=3\), it is selected in 46 of the 50 runs
(\(92\%\), \(95\%\) Wilson interval \([81,97]\%\)), and~the mean \(\Delta\mathrm{RMSE}\),
\(3.45\times10^{-3}\) with \(95\%\) confidence interval \mbox{\([2.86,4.03]\times10^{-3}\),} is
clearly positive. In~Case B with \(\alpha=1\), it is selected in 30 runs (\(60\%\),
\([46,72]\%\)), with~10 losses and 10 ties; the wins outnumber the losses, but~the
confidence interval of the mean difference, \([-0.01,0.79]\times10^{-3}\), includes zero,
so a mean advantage is not established. In~this setting, 10 of the 50 fits reach the
boundary \(\alpha_{\mathrm{fit}}=1\). The~supporting sign and Wilcoxon tests for all noisy
settings are reported in the Supplementary~Material.

Table~\ref{tab:noisy_three_family} reports the three-family comparison. Under~noise, the~biexponential no longer dominates the classical Case C as it did without noise: the
Mittag--Leffler law is selected in 30 of the 50 runs (\(60\%\), \([46,72]\%\)), the~biexponential in 17, and~the exponential in 3, and~the Mittag--Leffler model has, on~average, a~lower held-out RMSE than the biexponential. In~Case B with \(\alpha=1\), the~selections are spread across the three families, with~7 three-way ties in which the fitted
families coincide numerically. In~the fractional reference, Case B with \(\alpha=0.8\), the~Mittag--Leffler law is selected in 36 of the 50 runs (\(72\%\), \([58,83]\%\)). The~Mittag--Leffler selection frequencies of the classical Case C and of the fractional
reference therefore have overlapping intervals. In~all three settings, the~biexponential
usually attains the lowest unpenalized fitting-horizon RMSE, whereas the
horizon-extrapolation rankings often~differ.
\vspace{-20pt}

\vspace{-6pt}

\begin{table}[H]
	\caption{Three-family comparison under observational noise (\(\sigma_j=0.03|x_{0,j}|\),
		fifty seeds). Counts of runs in which each family has the lowest held-out RMSE, with~		\(95\%\) Wilson intervals for the selection frequency; ties: lowest RMSE shared within
		\(10^{-8}\). \(\Delta_{\mathrm{ML,Bi}}=\mathrm{RMSE}^{*}_{\mathrm{ML}}-
		\mathrm{RMSE}^{*}_{\mathrm{Bi}}\), mean with \(95\%\) confidence interval, in~units of
		\(10^{-3}\); negative values favor the Mittag--Leffler model.}
	\label{tab:noisy_three_family}

\begin{adjustwidth}{0pt}{0pt}
\begin{tabularx}{\fulllength}{lcCccC}
		\toprule
		\textbf{Setting} &\textbf{ Exp} & \textbf{ML} & \textbf{Bi} & \textbf{Ties} & \boldmath{\(\Delta_{\mathrm{ML,Bi}}\)} \\
		\midrule
		Classical, Case B, \(\alpha=1\) & 8 [0.08, 0.29] & 15 [0.19, 0.44] & 20 [0.28, 0.54] & 7
		& \(-0.14\) [\(-0.36\), 0.08] \\
		Classical, Case C, \(s=3\) & 3 [0.02, 0.16] & 30 [0.46, 0.72] & 17 [0.22, 0.48] & 0
		& \(-0.69\) [\(-1.17\), \(-0.22\)] \\
		Fractional, Case B, \(\alpha=0.8\) & 1 [0.00, 0.10] & 36 [0.58, 0.83] & 13 [0.16, 0.40] & 0
		& \(-1.56\) [\(-2.16\), \(-0.97\)] \\
		\bottomrule
	\end{tabularx}
\end{adjustwidth}
\end{table}
\unskip

\subsection{Stable~Sensitivity}
\label{Sec:SensitivityResults}

To test whether the preference observed for heterogeneous Case B depends on one
particular analysis configuration, we examine the Case B sensitivity study at \(\alpha=0.8\).  It shows no
change in the qualitative conclusion across the tested variations of \(K\), history length,
fitting interval, and~sampling interval.  Across all configurations,
\(\alpha_{\mathrm{fit}}\) remains between \(0.7756\) and \(0.7823\), and~the
Mittag--Leffler fit has lower horizon-extrapolation RMSE than the single
exponential in all ten seeds.  The~complete numerical results are reported in
Supplementary Table~S7.

As specified in Section~\ref{Sec:SensitivityProtocol}, the~sampling-interval rows
preserve the physical windows rather than \(m_{\mathrm{hist}}\), and~\(\Delta\mathrm{AIC}\) is compared only between configurations with the same number of
fitted~points.

\subsection{Numerical Reliability of the Single-Record~Example}
\label{Sec:NumericalReliability}

In its original analysis, the~single-record example was screened through its kNN
forecast-error growth profile. Over~the pre-transition interval, a~saturating logistic
curve, a~classical sigmoidal reference whose relation to Mittag--Leffler-based sigmoids has
been studied in~\cite{Hristov2025Logistic,Hristov2026Sigmoids}, described that profile
better than the exponential and Mittag--Leffler laws, but~the ranking changed with the
solver step (Supplementary Material). The~verification of
Section~\ref{Sec:NumericalCheck} shows~why.

\textls[-15]{The predictor-only mode of \texttt{fde12} reproduces the original explicit runs: at
\mbox{\(\Delta t=0.01\),}} \(0.005\), and~\(0.0025\), the~post-transient standard deviations agree
with those of the original computations to within \(10^{-4}\). As~the step is refined,
however, this regime changes qualitatively (Table~\ref{tab:numerical_reliability}): its
late-window standard deviation decreases from \(2.37\) at \(\Delta t=0.01\) to \(0.44\) at
\(\Delta t=0.00125\), and~its final state approaches the nonzero equilibrium \(E_-\). The~PECE solutions show substantially greater stability under refinement: their summary
diagnostics vary little over the same step-size range, and~the root-mean-square
differences between successive PECE solutions, \(0.034\), \(0.0084\), and~\(0.0022\),
decrease by factors close to four under step halving. The~smaller late-window variability
is consistent with a damped regime over the simulated time interval. An~iterated corrector
at \(\Delta t=0.005\) gives closely agreeing summary statistics, whereas the control at
\(\alpha=0.93\) retains a sustained large-amplitude oscillatory regime over the simulated
interval, with~a post-transient standard deviation of about \(2.44\); the
predictor--corrector scheme therefore does not generically suppress~oscillations.

\begin{table}[H]
	\caption{Predictor-only and PECE integration of System~\eqref{Eq:sene_system} at
		\(\alpha=0.916\). SD: standard deviation of \(x\) over \(t\geq40\) (post-transient) and
		over \(180\leq t\leq200\) (late window); distance: final distance from the nearest
		nonzero equilibrium.}
	\label{tab:numerical_reliability}
	 
	\begin{tabularx}{\textwidth}{lcCCcCC}
		\toprule
		& \multicolumn{3}{c}{\textbf{Predictor-Only}} & \multicolumn{3}{c}{\textbf{PECE}} \\
		\cmidrule{2-7}
		\boldmath{\multirow{-2}{*}{\(\Delta t\)}} & \textbf{SD, Post} & \textbf{SD, Late} & \textbf{Distance} & \textbf{SD, Post} & \textbf{SD, Late} & \textbf{Distance} \\
		\midrule
		0.01    & 2.345 & 2.367 & 5.04 & 0.929 & 0.276 & 0.492 \\
		0.005   & 2.353 & 2.322 & 3.80 & 0.926 & 0.274 & 0.489 \\
		0.0025  & 1.417 & 0.744 & 1.28 & 0.925 & 0.274 & 0.488 \\
		0.00125 & 1.117 & 0.444 & 0.577 & 0.925 & 0.273 & 0.488 \\
		\bottomrule
	\end{tabularx}
\end{table}

The local stability analysis is consistent with these observations. At~\(E_\pm\), the~Jacobian eigenvalues are approximately \(-3.482\) and \(0.241\pm1.902\mathrm{i}\), giving the
critical order \(\alpha_c=(2/\pi)\,|\arg(0.241+1.902\mathrm{i})|\approx0.91985\). The~order
\(\alpha=0.916\) lies just below this threshold, so the nonzero equilibria are locally
asymptotically stable, whereas \(\alpha=0.93\) lies above it. This local result does not
determine the global dynamics, but~together with the refinement evidence it shows that the
sustained oscillatory regime of the original record is not numerically robust. The~forecast-error profile and model ranking computed from that record are therefore withdrawn
as evidence about the continuous fractional system; they are retained in the Supplementary
Material only to document the sensitivity of the original~discretization.

\section{Discussion}
\label{Sec:Discussion}
\unskip

\subsection{Controlled Stable~Experiments}
\label{Sec:ControlledStableDiscussion}
The comparison of profile constructions (Section~\ref{Sec:BaselineResults}) clarifies
the role of kNN in the stable screening results. In~all tested noise-free
configurations, the~kNN, constant training-mean, and~direct-average profiles led
to the same exponential-versus--Mittag--Leffler decision, although~their fitted
orders differed slightly. This agreement suggests that features of the aggregate
relaxation profiles sufficient for the binary screening decision are preserved
across the three constructions under these conditions; the present experiments
do not identify which features are responsible. The~kNN construction is therefore
not necessary to reproduce these screening decisions, but~it remains a useful
practical construction: among the three alternatives tested, it alone is both
out of sample and conditioned on the observed history of each test trajectory.
It can therefore be applied to new trajectories when their histories, but~not
their futures, are available, and~it predicts individual trajectories far more
accurately than the constant predictor. Part of that advantage reflects the
amplitude information contained in each history, and~the comparison does not
separate it from information about relaxation~shape.

The near-exact recovery obtained in the matched Case A calibration
(Section~\ref{Sec:ExactCalibration}) is expected under the conditions of that
experiment: the generating and fitted families
coincide, the~trajectories are noise-free, and~the ensemble is
multiplicatively separable. The~result therefore provides a calibration check of the
numerical implementation under matched conditions, but~provides little evidence for
general fractional-order~identifiability.

Case B (Section~\ref{Sec:HeterogeneousRates}) changes the interpretation
because the ensemble profile is no longer generated by a single
Mittag--Leffler law: individual trajectories share a
common fractional order but carry different relaxation rates. The~fitted
ensemble-level order is systematically shifted downward relative to the
trajectory-level generating order: the heterogeneous mixture yields an ensemble-level
fitted shape parameter below the common generating order of its components. This
offset should not be read as an estimation error relative to a uniquely defined
ensemble-level order, because~no such single-order generating law exists for
the mixed~profile.

The same direction of departure becomes more pronounced in the classical
distributed-rate control of Case C (Section~\ref{Sec:SpecificityResults}), where rate heterogeneity alone drives the
fitted shape parameter well below the integer-order limit even though every
trajectory is purely exponential. This progression motivates the specificity
analysis in Section~\ref{Sec:SpecificityDiscussion}, where the role of the
candidate model set is considered~explicitly.

\subsection{Specificity and Candidate-Set~Breadth}
\label{Sec:SpecificityDiscussion}

Case C provides the direct test of candidate-set dependence.  In~this classical
distributed-rate control (Section~\ref{Sec:SpecificityResults}), the~initial exponential-versus-Mittag--Leffler comparison favors the Mittag--Leffler family even though every
trajectory is generated by a classical exponential mechanism.  Without~noise, that
preference disappears once the biexponential family is included.  The~binary
comparison therefore produced a candidate-set-dependent positive screen: the
narrower candidate set made a classical mixture appear
Mittag--Leffler-compatible.  For~\(0<\alpha<1\), \(E_\alpha(-t^\alpha)\) is
completely monotone and admits a nonnegative spectral representation as a
continuous superposition of exponential relaxations
~\cite{Mainardi2000On,Mainardi2020Why}.  Relaxation shape alone is therefore, in~general, insufficient to identify fractional origin uniquely; this provides a
mathematical counterpart to the candidate-set dependence observed in Case~C.

Case C thus directly demonstrates that model preference can change with candidate-set
breadth.  The~monotone decrease of
\(\alpha_{\mathrm{fit}}\) with increasing rate spread in Case C is consistent
with the fitted shape parameter absorbing unresolved rate heterogeneity, but~the present experiments do not establish that interpretation as a general
mechanism.

Under observational noise, this recovery of specificity is incomplete
(Section~\ref{Sec:NoiseResults}). With~the biexponential admitted, the~Mittag--Leffler law is still selected in 30 of 50 runs for the classical
Case C ensemble, and~it has, on~average, the~lower held-out error,
much as in the fractional reference, where it is selected in 36 of
50 runs. Although~these frequencies do not establish equivalent
selection probabilities, their proximity in the tested configurations
indicates that Mittag--Leffler selection alone provides limited
discrimination between classical rate heterogeneity and fractional
relaxation. The~biexponential attains the lowest fitting-horizon error
in most runs yet is selected less often on the held-out horizons;
this contrast is consistent with part of its additional flexibility
being spent on the noise, although~the experiments do not isolate that
mechanism. Broadening the candidate set therefore reduces
false-positive screens but does not, by~itself, guarantee discrimination
once the observations are noisy; the results indicate limited
specificity under the tested conditions, not a general
indistinguishability of the two~mechanisms.

Operationally, a~positive Mittag--Leffler screen should therefore
be challenged, where possible, against~flexible nonfractional alternatives
before being interpreted as evidence favoring a fractional generating
mechanism, although~under observational noise such a challenge may not
suffice.  No finite candidate set can be exhaustive, so even after such
comparisons a preference supports compatibility within the tested model set
rather than unique mechanism~identification.

\subsection{Numerical Provenance of Simulated~Data}
\label{Sec:NumericalProvenance}

Data-driven screening methods are often developed and tested on simulated observations,
and their conclusions can then be no more reliable than the simulations themselves. The~single-record example shows that numerical artifacts in simulated observations can
propagate into screening conclusions: a coarse explicit discretization produced a
sustained oscillatory regime on which the screening analysis returned a definite, though~resolution-dependent, model ranking, whereas step refinement and an independent
predictor--corrector scheme did not confirm that regime
(Section~\ref{Sec:NumericalReliability}). The~example is particularly instructive because the
chosen order, \(\alpha=0.916\), lies close to the local stability threshold
\(\alpha_c\approx0.91985\) of the nonzero equilibria, where small numerical errors can
plausibly have qualitative consequences. Simulated records intended for screening should
therefore be checked under step refinement and, where possible, with~an independent
numerical scheme, especially when the parameters lie close to a stability boundary.
Agreement across such checks does not by itself establish the scientific interpretation
of a screening outcome, whereas disagreement prevents confident attribution to the
continuous model without further numerical~verification.

\subsection{Scope and~Limitations}
\label{Sec:Limitations}

The proposed workflow is a screening procedure for finite-window profile
compatibility, not a mechanism-identification method: a positive screen indicates
compatibility relative to the tested alternatives, and~the fitted shape parameter should
not be interpreted as a physical fractional order unless that interpretation is supported
by independent modeling~evidence.

The nonfractional controls include heterogeneous classical rate mixtures,
biexponential alternatives, and~separable curved shape controls. They are
nevertheless not exhaustive. In~particular, nonlinear classical generators
and explicitly nonstationary alternatives remain outside the present study and
constitute natural extensions of the specificity analysis. More generally, no practical
analysis can include every competing family, so a positive screen cannot exclude all
nonfractional~explanations.

Observational noise introduces a further practical limitation.  The~Mittag--Leffler family converges to the exponential family as
\(\alpha\to1\), so discrimination necessarily vanishes at the integer-order
limit even for exact data.  Observational noise makes the residual difference
between the two families harder to resolve already before that limit is
reached, despite the fitted order remaining close to the generating value over
the fractional cases.  A~weak or vanishing preference near the integer-order
regime should therefore not be interpreted as evidence against either~family.

Finally, the~uncertainty analysis in the stable noisy experiment has two
different levels. Variation across repeated realizations is summarized over
independent seeds, whereas the trajectory bootstrap for the representative run
is conditional on the estimated training predictor and fixed training set.
The latter is therefore not a full-pipeline uncertainty~estimate.

The analyses added in this study have their own scope. The~comparison of profile
constructions covers thirteen noise-free configurations and includes no amplitude-aware
predictor, so it neither extends automatically to noisy profiles nor separates amplitude
from shape information. The~noisy specificity analysis uses a single observational-noise
level and three targeted generating settings, and~its fractional reference, Case B with
\(\alpha=0.8\), does not establish general sensitivity. Several noise-free analyses,
including the comparison of profile constructions and the stable sensitivity study, rest on
ten seeds per configuration. Formal paired hypothesis tests for fitted-order differences
have not been conducted; the corresponding estimates and intervals are interpreted
descriptively, particularly near the constrained integer-order boundary. Finally, the~numerical-reliability check concerns one
nonlinear system, one initial condition, and~a limited selection of fractional orders.
These limits bound the strength of any inferential or mechanistic claim drawn from the
fitted profiles and should be taken into account when the method is used as a comparative
screening~tool.

\section{Conclusions}
\label{Sec:Conclusions}

This study develops and evaluates an equation-free, finite-window compatibility screen
based on multi-horizon relaxation profiles constructed by kNN prediction and on the
comparative fitting of reference temporal laws. The~results support a deliberately limited
interpretation of that screen. In~the matched, noise-free calibration, the~fitted
Mittag--Leffler family reproduces the generating behavior to numerical precision; this
result is a calibration check rather than evidence of general fractional-order
identifiability.

When the individual trajectories retain a common fractional order but differ only in their
relaxation rates, the~ensemble-level fitted shape parameter moves systematically below the
common generating order. The~same directional effect becomes stronger for purely classical
distributed-rate ensembles. The~present results therefore show that relaxation-rate
heterogeneity alone can produce an apparently Mittag--Leffler-compatible ensemble profile
even when the individual components are purely~exponential.

The central result is that model preference is conditional on the candidate set and on the
observation conditions. Without~noise, the~apparent Mittag--Leffler preference obtained for
a classical distributed-rate ensemble against a single exponential disappears when a
flexible biexponential alternative is added. Under~observational noise, broadening the
candidate set reduces such false-positive screens but does not remove them: the
Mittag--Leffler law is still selected in a majority of runs for the classical ensemble, at~a
frequency close to that obtained for a fractional reference. A~positive Mittag--Leffler
screen should therefore be interpreted only relative to the candidate families actually
tested and the conditions of observation, and~not as unique identification of the mechanism
generating the~data.

In the tested noise-free configurations, the~binary screening decision is reproduced by
simpler profile constructions, including a constant training-mean predictor and the
directly observed ensemble average. The~kNN construction is therefore not necessary for
these decisions; its value lies in providing an out-of-sample, history-conditioned profile
and substantially more accurate individual-trajectory predictions than the constant
training-mean~predictor.

\clearpage
The single-record example, initially investigated as a chaotic test case, is retained only
as a cautionary illustration. Step refinement and an independent predictor--corrector
scheme showed that the sustained regime produced by the original explicit discretization is
not numerically robust, and~its screening results are withdrawn as evidence about the
continuous system. Screening conclusions drawn from simulated data therefore require
careful verification of the numerical procedure used to generate those~data.

Taken together, the~results support the use of the method as a comparative finite-window
compatibility screen. Its outcomes must be interpreted relative to the candidate set, with~attention to relaxation-rate heterogeneity, and~kept distinct from physical fractional-order
identification; observational noise weakens discrimination both near the integer-order
limit and between classical mixtures and fractional~relaxation.

\vspace{6pt}

\section*{Supplementary Material}
The following supporting information can be
	downloaded at: \linksupplementary{s1}.
	Table S1: Per-start biexponential optimizer diagnostics for
	Cases A and B at $\alpha=0.8$ and Case C at $s=3$ (main run);
	Table S2: Fitted parameters for the expanded candidate comparison
	of Table 5 (mean $\pm$ standard deviation over ten seeds);
	Table S3: Comparison of profile constructions in all thirteen
	noise-free configurations (ten seeds each);
	Table S4: Binary exponential-versus--Mittag--Leffler comparison
	under observational noise;
	Table S5: Paired held-out RMSE differences in the noisy
	three-family comparison;
	Table S6: Biexponential optimizer diagnostics in the noisy
	three-family comparison;
	Table S7: One-at-a-time sensitivity analysis on Case B
	at $\alpha=0.8$;
	Table S8: Reproduction of the original explicit runs by the
	predictor-only mode of \texttt{fde12};
	Table S9: Predictor--corrector (PECE) refinement at $\alpha=0.916$;
	Table S10: Optimization settings of the original single-record fits;
	Table S11: Physical-time split of the original record;
	Table S12: Pre-transition model comparison for the original
	kNN forecast-error growth profile using the physical-time split
	with embargo;
	Table S13: Optimizer diagnostics for the original pre-transition
	fits at the three solver steps;
	Table S14: Sensitivity of the original pre-transition comparison
	to the training allocation of the physical-time split;
	Table S15: Solver-step comparison for the original fractional system;
	Table S16: Original record, $\Delta t=0.01$: sensitivity of the
	pre-transition fits to the fitting window and to the error metric;
	Figure S1: kNN forecast-error growth profile of the original
	record (baseline split) with the exponential, Mittag--Leffler,
	and logistic laws fitted over the pre-transition interval;
	Figure S2: Identifiability of the Mittag--Leffler order for the
	baseline original fit;
	Figure S3: Solver-step comparison for the original fractional system.
	The Supplementary Material also cites the work of
	Apelblat~\cite{Apelblat2010Integral}.

\section*{Author Contributions}
Conceptualization, N.N., A.V. and B.C.; methodology, N.N., A.V. and B.C.;
	software, N.N. and A.V.; validation, N.N., A.V. and B.C.; formal analysis, N.N., A.V. and B.C.;
	investigation, N.N. and A.V.; data curation, N.N., A.V. and B.C.; visualization, N.N.;
	writing---original draft preparation, N.N. and A.V.; writing---review and editing, N.N., A.V. and B.C.; supervision, B.C.; funding acquisition, B.C.
	All authors have read and agreed to the published version of the~manuscript.

\section*{Funding}
\cbl{This research was supported by the European Regional Development Fund (ERDF) 2021--2027 under Project AI4AM---EFRE1052. This work was supported by the Open Access Publishing Fund of the Free University of~Bozen-Bolzano.}

\section*{Institutional Review Board Statement}
Not applicable.

\section*{Informed Consent Statement}
Not applicable.

\section*{Data Availability Statement}
The numerical experiments in this study are fully synthetic. The~MATLAB scripts, generated data, and~optimizer logs required to reproduce the reported results are openly available on Zenodo at \url{https://doi.org/10.5281/zenodo.23494805}. The~repository README files map each script and output file to the corresponding tables and figures of the manuscript and the Supplementary Material. The~Mittag--Leffler function is evaluated with the MATLAB implementation of Garrappa~\cite{Garrappa2015Numerical}, which is not redistributed; the README indicates where to obtain it.

\section*{Acknowledgments}
\cbl{Bruno Carpentieri is a member of the
		Gruppo Nazionale per il Calcolo Scientifico (GNCS)
		of the Istituto Nazionale di Alta Matematica (INdAM).}

\section*{Conflicts of Interest}
The authors declare no financial conflicts of interest. B.C. serves as a Guest Editor of a Special Issue in \emph{Fractal and~Fractional}; he will not participate in the editorial handling, peer review, or~decision-making process for this manuscript. The~funders had no role in the design of the study; in the collection, analyses, or~interpretation of data; in the writing of the manuscript; or in the decision to publish the~results.

\section*{Abbreviations}
~ \vspace{-4pt}

 \noindent 
\begin{tabular}{@{}ll}

 AIC & Akaike information criterion\\
BIC & Bayesian information criterion\\
FDE & fractional differential equation\\
kNN & k-nearest neighbors\\
ML & Mittag--Leffler\\
RMSE & root mean square error
\end{tabular}


\clearpage
\begin{thebibliography}{999}

\bibitem[Podlubny(1999)]{podlubny1999fractional}
Podlubny, I.
\newblock {\em Fractional Differential Equations}; Mathematics
  in Science and Engineering; Academic Press: San Diego,  CA, USA, 1999; Volume 198. 

\bibitem[Diethelm(2010)]{diethelm2010analysis}
Diethelm, K.
\newblock {\em The Analysis of Fractional Differential Equations}; 
  Lecture Notes in Mathematics; Springer: Berlin,  2010; Volume 2004.
\newblock {\url{https://doi.org/10.1007/978-3-642-14574-2}}.

\bibitem[Herrmann(2014)]{Herrmann2014Fractional}
Herrmann, R.
\newblock {\em Fractional Calculus: An Introduction for Physicists}, 2nd ed.;
  World Scientific: Singapore,  2014.
\newblock {\url{https://doi.org/10.1142/8934}}.

\bibitem[Mainardi and Gorenflo(2000)]{Mainardi2000On}
Mainardi, F.; Gorenflo, R.
\newblock On {Mittag-Leffler}-type functions in fractional evolution processes.
\newblock {\em J. Comput. Appl. Math.} {\bf 2000},
  {\em 118},~283--299.
\newblock {\url{https://doi.org/10.1016/S0377-0427(00)00294-6}}.

\bibitem[Mainardi(2020)]{Mainardi2020Why}
Mainardi, F.
\newblock Why the {Mittag-Leffler} function can be considered the queen
  function of the fractional calculus?
\newblock {\em Entropy} {\bf 2020}, {\em 22},~1359.
\newblock {\url{https://doi.org/10.3390/e22121359}}.

\bibitem[Gorenflo et~al.(2020)Gorenflo, Kilbas, Mainardi, and
  Rogosin]{Gorenflo2020MittagLeffler}
Gorenflo, R.; Kilbas, A.A.; Mainardi, F.; Rogosin, S.
\newblock {\em {Mittag-Leffler} Functions, Related Topics and Applications},
  2nd ed.; Springer: Berlin/Heidelberg, Germany, 
  2020.
\newblock {\url{https://doi.org/10.1007/978-3-662-61550-8}}.

\bibitem[Metzler and Klafter(2002)]{Metzler2002From}
Metzler, R.; Klafter, J.
\newblock From stretched exponential to inverse power-law: Fractional dynamics,
  {Cole-Cole} relaxation processes, and beyond.
\newblock {\em J. Non-Cryst. Solids} {\bf 2002}, {\em
  305},~81--87.
\newblock {\url{https://doi.org/10.1016/S0022-3093(02)01124-9}}.

\bibitem[Granger and Joyeux(1980)]{Granger1980AN}
Granger, C.W.J.; Joyeux, R.
\newblock An introduction to long-memory time series models and fractional
  differencing.
\newblock {\em J. Time Ser. Anal.} {\bf 1980}, {\em 1},~15--29.
\newblock {\url{https://doi.org/10.1111/j.1467-9892.1980.tb00297.x}}.

\bibitem[Delle~Monache et~al.(2015)Delle~Monache, Grassi, and Santucci~de
  Magistris]{Monache2015Testing}
Delle~Monache, D.; Grassi, S.; Santucci~de Magistris, P.
\newblock \textit{Testing for Level Shifts in Fractionally Integrated Processes: A
  State Space Approach};
\newblock CREATES Research Paper 2015-30; Department of Economics and Business
  Economics, Aarhus University: Aarhus,  Denmark, 2015.

\bibitem[Mainardi(2018)]{Mainardi2017A}
Mainardi, F.
\newblock A note on the equivalence of fractional relaxation equations to
  differential equations with varying coefficients.
\newblock {\em Mathematics} {\bf 2018}, {\em 6},~8.
\newblock {\url{https://doi.org/10.3390/math6010008}}.

\bibitem[Poinot and Trigeassou(2004)]{Poinot2004Identification}
Poinot, T.; Trigeassou, J.C.
\newblock Identification of fractional systems using an output-error technique.
\newblock {\em Nonlinear Dyn.} {\bf 2004}, {\em 38},~133--154.
\newblock {\url{https://doi.org/10.1007/s11071-004-3751-y}}.

\bibitem[Li and Rosenfeld(2022)]{Li2021Fractional}
Li, X.; Rosenfeld, J.A.
\newblock Fractional order system identification with occupation kernel
  regression.
\newblock {\em IEEE Control Syst. Lett.} {\bf 2022}, {\em 6},~19--24.
\newblock {\url{https://doi.org/10.1109/LCSYS.2020.3046408}}.

\bibitem[Zhang et~al.(2024)Zhang, Lu, Liu, and Liu]{Zhang2024Sparse}
Zhang, T.; Lu, Z.R.; Liu, J.K.; Liu, G.
\newblock Sparse identification of fractional chaotic systems based on the
  time-domain data.
\newblock {\em Chin. J. Phys.} {\bf 2024}, {\em 89},~160--173.
\newblock {\url{https://doi.org/10.1016/j.cjph.2024.02.050}}.

\bibitem[Duh{\'e} et~al.(2024)Duh{\'e}, Victor, Melchior, Abdelmoumen, and
  Roubertie]{duhe2024recursive}
Duh{\'e}, J.F.; Victor, S.; Melchior, P.; Abdelmoumen, Y.; Roubertie, F.
\newblock Recursive system identification of continuous-time fractional systems
  for all parameter estimation.
\newblock {\em IEEE Trans. Control Syst. Technol.} {\bf 2024},
  {\em 32},~2037--2049.
\newblock {\url{https://doi.org/10.1109/TCST.2024.3407580}}.

\bibitem[Yuan et~al.(2024)Yuan, Xu, Liu, Wang, and Lu]{Yuan2024Identification}
Yuan, M.; Xu, W.; Liu, F.; Wang, L.; Lu, Y.
\newblock Identification method for a fractional-order system in terms of
  equivalent dynamic properties.
\newblock {\em Chaos} {\bf 2024}, {\em 34},~073119.
\newblock {\url{https://doi.org/10.1063/5.0187031}}.

\bibitem[Ben~Brahim et~al.(2024)Ben~Brahim, El~Alaoui, and
  Zguaid]{BenBrahim2024Some}
Ben~Brahim, H.; El~Alaoui, F.Z.; Zguaid, K.
\newblock Some results regarding observability and initial state reconstruction
  for time-fractional systems.
\newblock {\em  Int. J. Optim. Control. Theor. Appl. (IJOCTA)} {\bf 2024}, {\em 14},~99--112.
\newblock {\url{https://doi.org/10.11121/ijocta.1468}}.

\bibitem[Du et~al.(2013)Du, Wang, and Hu]{Du2013Measuring}
Du, M.; Wang, Z.; Hu, H.
\newblock Measuring memory with the order of fractional derivative.
\newblock {\em Sci. Rep.} {\bf 2013}, {\em 3},~3431.
\newblock {\url{https://doi.org/10.1038/srep03431}}.

\bibitem[Tarasov(2019)]{Tarasov2019On}
Tarasov, V.E.
\newblock On history of mathematical economics: Application of fractional
  calculus.
\newblock {\em Mathematics} {\bf 2019}, {\em 7},~509.
\newblock {\url{https://doi.org/10.3390/math7060509}}.

\bibitem[Wei and Yuan(2024)]{Wei2024Identifying}
Wei, Z.; Yuan, N.
\newblock Identifying climate memory impacts on climate network analysis: A new
  approach based on fractional integral techniques.
\newblock {\em Clim. Dyn.} {\bf 2024}, {\em 62},~3465--3476.
\newblock {\url{https://doi.org/10.1007/s00382-023-07076-z}}.

\bibitem[Velichko et~al.(2025)Velichko, Belyaev, and Boriskov]{Velichko2025A}
Velichko, A.; Belyaev, M.; Boriskov, P.
\newblock A novel approach for estimating largest {Lyapunov} exponents in
  one-dimensional chaotic time series using machine learning.
\newblock {\em Chaos} {\bf 2025}, {\em 35},~101101.
\newblock {\url{https://doi.org/10.1063/5.0289352}}.

\bibitem[Garrappa(2015)]{Garrappa2015Numerical}
Garrappa, R.
\newblock Numerical evaluation of two and three parameter {Mittag-Leffler}
  functions.
\newblock {\em SIAM J. Numer. Anal.} {\bf 2015}, {\em
  53},~1350--1369.
\newblock {\url{https://doi.org/10.1137/140971191}}.

\bibitem[Sene(2021)]{sene2021analysis}
Sene, N.
\newblock Analysis of a fractional-order chaotic system in the context of the
  {Caputo} fractional derivative via bifurcation and {Lyapunov} exponents.
\newblock {\em J. King Saud Univ.-Sci.} {\bf 2021}, {\em
  33},~101275.
\newblock {\url{https://doi.org/10.1016/j.jksus.2020.101275}}.

\bibitem[Sene(2022)]{Sene2022Theory}
Sene, N.
\newblock Theory and applications of new fractional-order chaotic system under
  {Caputo} operator.
\newblock {\em  Int. J. Optim. Control. Theor. Appl. (IJOCTA)} {\bf 2022}, {\em 12},~20--38.
\newblock {\url{https://doi.org/10.11121/ijocta.2022.1108}}.

\bibitem[Diethelm et~al.(2005)Diethelm, Ford, Freed, and
  Luchko]{Diethelm2005Algorithms}
Diethelm, K.; Ford, N.J.; Freed, A.D.; Luchko, Y.
\newblock Algorithms for the fractional calculus: A selection of numerical
  methods.
\newblock {\em Comput. Methods Appl. Mech. Eng.} {\bf
  2005}, {\em 194},~743--773.
\newblock {\url{https://doi.org/10.1016/j.cma.2004.06.006}}.

\bibitem[Hristov(2025)]{Hristov2025Logistic}
Hristov, J.
\newblock The logistic models and sigmoid functions: A variety of models and
  solution perspectives.
\newblock {\em Trans. Comput. Model. Intell. Syst.}
  {\bf 2025}, {\em 1},~10002.
\newblock {\url{https://doi.org/10.65112/tcmis.10002}}.

\bibitem[Hristov(2026)]{Hristov2026Sigmoids}
Hristov, J.
\newblock Sigmoids based on {Mittag-Leffler} functions: Ideas, modeling, and
  computational experiments.
\newblock {\em Trans. Comput. Model. Intell. Syst.}
  {\bf 2026}, {\em 3},~10073.
\newblock {\url{https://doi.org/10.65112/tcmis.10073}}.

\bibitem[Apelblat(2010)]{Apelblat2010Integral}
Apelblat, A.
\newblock {\em Integral Transforms and Volterra Functions};
Nova Science Publishers: New York, NY, USA, 2010.

\end{thebibliography}
\end{document}